\documentclass[american,a4paper,10pt,leqno]{amsart}

\usepackage{babel,enumerate,amssymb,amsthm}
\usepackage[arrow,curve,matrix,tips,2cell]{xy}
  \SelectTips{eu}{10} \UseTips
  \UseAllTwocells

\DeclareMathAlphabet{\matheurm}{U}{eur}{m}{n}
\newcommand{\cell}{\operatorname{cell}}

\newcommand{\iso}{\operatorname{iso}}
\newcommand{\op}{\operatorname{op}}
\newcommand{\pt}{\operatorname{pt}}
\newcommand{\per}{\operatorname{per}}
\newcommand{\sing}{\operatorname{sing}}

\newcommand{\Jdtr}{\mathbf{dtr}}
\newcommand{\Jntr}{\mathbf{ntr}}

\newcommand{\Ab}{\matheurm{Ab}}

\newcommand{\Cat}{\matheurm{Cat}}
\newcommand{\FGINJ}{\matheurm{FGINJ}}
\newcommand{\Groupoids}{{\matheurm{Groupoids}}}
\newcommand{\Or}{\matheurm{Or}\hspace*{.1em}}   %% I HAVE MODIFIED THIS MACRO %%

\newcommand{\Sp}{\matheurm{Sp}}

\DeclareMathOperator{\Aut}{Aut}
\DeclareMathOperator{\aut}{aut}
\DeclareMathOperator{\ch}{ch}
\DeclareMathOperator{\coequ}{coequ}
\DeclareMathOperator{\coker}{coker}
\DeclareMathOperator{\con}{con}
\DeclareMathOperator{\diag}{diag}
\DeclareMathOperator{\dtr}{dtr}
\DeclareMathOperator{\ntr}{ntr}
\DeclareMathOperator{\id}{id}
\DeclareMathOperator{\ind}{ind}
\DeclareMathOperator{\Idem}{Idem}
\DeclareMathOperator{\map}{map}
\DeclareMathOperator{\mor}{mor}
\DeclareMathOperator{\obj}{obj}
\DeclareMathOperator{\pr}{pr}
\DeclareMathOperator{\res}{res}
\DeclareMathOperator{\sub}{sub}
\DeclareMathOperator{\tr}{tr}
\DeclareMathOperator{\Seg}{Seg}
\DeclareMathOperator{\Sw}{Sw}

\DeclareMathOperator{\Tot}{Tot}

\newcommand{\negspace}{\hspace*{-.03em}} %% THIS IS A NEW MACRO %%
\DeclareMathOperator{\HH}{H\negspace H}  %% THIS IS A NEW MACRO %%
\DeclareMathOperator{\HC}{H\negspace C}  %% THIS IS A NEW MACRO %%
\DeclareMathOperator{\HN}{H\negspace N}  %% THIS IS A NEW MACRO %%
\DeclareMathOperator{\HP}{H\negspace P}  %% THIS IS A NEW MACRO %%
\DeclareMathOperator{\HX}{H\negspace X}  %% THIS IS A NEW MACRO %%
\DeclareMathOperator{\CN}{C\negspace N}  %% THIS MACRO HAS BEEN MODIFIED %%
\DeclareMathOperator{\DK}{D\negspace\negspace K}  %% THIS MACRO HAS BEEN MODIFIED %%

\DeclareMathOperator*{\colim}{colim}
\DeclareMathOperator*{\sma}{\wedge}  % \smash already defined

\DeclareMathOperator*{\tensor}{\otimes}

\newcommand{\Fin}  {\mathcal{F}\text{in}}
\newcommand{\FCyc}{\mathcal{FC}\text{yc}}
\newcommand{\VCyc}{\mathcal{VC}\text{yc}}

  \newcommand{\IA}{\mathbb{A}}
  
  \newcommand{\IC}{\mathbb{C}}
  
  \newcommand{\IE}{\mathbb{E}}
  
  \newcommand{\IH}{\mathbb{H}}

  \newcommand{\IN}{\mathbb{N}}

  \newcommand{\IQ}{\mathbb{Q}}

  \newcommand{\IZ}{\mathbb{Z}}

  \newcommand{\JC}{\mathbf{C}}    
  \newcommand{\JD}{\mathbf{D}}    
  \newcommand{\JE}{\mathbf{E}}    
  \newcommand{\JF}{\mathbf{F}}    \newcommand{\Jf}{\mathbf{f}}
      
  \newcommand{\JH}{\mathbf{H}}    \newcommand{\Jh}{\mathbf{h}}

  \newcommand{\JK}{\mathbf{K}}

  \newcommand{\JN}{\mathbf{N}}    
      
  \newcommand{\JP}{\mathbf{P}}

  \newcommand{\JS}{\mathbf{S}}

  \newcommand{\JX}{\mathbf{X}}

  \newcommand{\cala}{\mathcal{A}}
  \newcommand{\calb}{\mathcal{B}}
  \newcommand{\calc}{\mathcal{C}}

  \newcommand{\calf}{\mathcal{F}}
  \newcommand{\calg}{\mathcal{G}}
  \newcommand{\calh}{\mathcal{H}}

  \newcommand{\caln}{\mathcal{N}}
  \newcommand{\calo}{\mathcal{O}}

\newcommand{\EGF}[2]{E_{#2}(#1)}               % klassifizierender Raum einer Familie CW-version

\newcommand{\conjclass}[1]{(#1)} % or: \newcommand{\conjclass}[1]{[#1]}
\DeclareMathOperator{\assembly}{assembly}
\newcommand{\into}{\hookrightarrow}

\DeclareMathOperator{\h}{h}
\newcommand{\minusinfinity}{{-\infty}}
\newcommand{\downdots}{\raisebox{-.55em}{$\ldots$}}
\newcommand{\stackover}[2]{\renewcommand{\arraystretch}{0.2}
\begin{array}[t]{c}
#2\\ {\scriptstyle #1}\\
\end{array}
\renewcommand{\arraystretch}{1}}

\newcounter{commentcounter}

\theoremstyle{plain}      \newtheorem{theorem}{Theorem}[section]
                          
                          \newtheorem{lemma}[theorem]{Lemma}
                          \newtheorem{corollary}[theorem]{Corollary}
                          \newtheorem{proposition}[theorem]{Proposition}
                          \newtheorem{addendum}[theorem]{Addendum}

\theoremstyle{definition} 
                          \newtheorem{example}[theorem]{Example}
                          \newtheorem{question}[theorem]{Question}
                          \newtheorem{remark}[theorem]{Remark}
                          \newtheorem{notation}[theorem]{Notation}

\makeatletter\let\c@equation=\c@theorem\makeatother

  \title[Detecting $K$-theory by cyclic homology]
              {Detecting $K$-Theory by cyclic homology}
       \author{Wolfgang L\"uck}
              \author{Holger Reich}
      \address{Westf\"alische Wilhelms-Universit\"at M\"unster\\
               Mathematisches Institut\\
               Einsteinstr.~62,
               D-48149 M\"unster, Germany}
        \email{lueck@math.uni-muenster.de}
      \urladdr{http://www.math.uni-muenster.de/u/lueck}
        \email{reichh@math.uni-muenster.de}
      \urladdr{http://www.math.uni-muenster.de/u/reichh}
               \date{\today}
     \keywords{}
    \subjclass[2000]{}

\begin{document}

\maketitle

\begin{abstract}
We discuss which part of the rationalized algebraic $K$-theory of a group ring is detected
via trace maps to Hochschild homology, cyclic homology, periodic cyclic or negative cyclic homology.

\smallskip
\noindent
Key words: algebraic K-theory of group rings, Hochschild homology, cyclic homology, trace maps.\\
Mathematics Subject Classification 2000: 19D55.
\end{abstract}

Dedicated to memory of Michel Matthey.

%%%%%%%%%%%%%%%%%%%%%%%%%%%%%%%%%%%%%%%%%%%%%%%%%%%%%%%%%%%%%%%%%%%%%%%%%%%%%%%%%%%%%%%%%%%%%%%%%%%%%%
%%%%%%%%%%%%%%%%%%%%%%%%%%%%%%%%%%%%%%%%%%%%%%%%%%%%%%%%%%%%%%%%%%%%%%%%%%%%%
\setcounter{section}{-1}
\section{Introduction and statement of results}

Fix a commutative ring $k$, referred to as the \emph{ground ring}. Let $R$ be a $k$-algebra, i.e.\ an associative ring $R$ together
with a unital ring homomorphism from $k$ to the center of $R$. We
denote by $\HH^{\otimes_k}_{\ast} ( R )$ the \emph{Hochschild homology} of $R$ relative to the ground ring $k$,
and similarly by $\HC^{\otimes_k}_{\ast}(R)$, $\HP^{\otimes_k}_{\ast}(R)$ and $\HN_{\ast}^{\otimes_k} (R)$
the \emph{cyclic}, the \emph{periodic cyclic} and the \emph{negative cyclic homology} of $R$ relative to $k$.
Hochschild homology receives a map from the algebraic $K$-theory,
which is known as the \emph{Dennis trace map}. There are variants of the Dennis trace taking values in cyclic,
periodic cyclic and negative cyclic homology (sometimes called \emph{Chern characters}), as displayed in the
following commutative diagram.
\begin{eqnarray} \label{trace-maps}
\vcenter{\xymatrix@C=3.2em{
& \HN_{\ast}^{\otimes_k} ( R ) \ar[r] \ar[d]^{\h} & \HP_{\ast}^{\otimes_k} ( R ) \ar[d] \\
K_{\ast} ( R ) \ar[ru]^-{\ntr} \ar[r]^-{\dtr} & \HH_{\ast}^{\otimes_k} ( R ) \ar[r] & \HC_{\ast}^{\otimes_k} ( R ).
         }}
\end{eqnarray}
For the definition of these maps, see \cite[Chapters~8~and~11]{Loday(1992)} and Section~\ref{sec:trace-maps} below.

In the following we will focus on the case of group rings $RG$, where $G$ is a group and we refer to the
$k$-algebra $R$ as the \emph{coefficient ring}.
We investigate the following question.
\begin{question} \label{que:main-question}
Which part of $K_{\ast}( R G ) \otimes_{\IZ} \IQ$ can be detected using linear trace invariants
like the Dennis trace to Hochschild homology, or its variants with values in cyclic homology,
periodic cyclic homology and negative cyclic homology\,?
\end{question}

For any group $G$, we prove ``\emph{detection results}'', which state that certain parts of $K_\ast(RG) \otimes_{\IZ} \IQ$
can be detected by the trace maps in diagram~\ref{trace-maps}, accompanied by ``\emph{vanishing results}'',
which state that a complement of the part which is then known to be detected is mapped to zero.
For the detection results, we only make assumptions on the coefficient ring $R$, whereas for
the vanishing results we additionally need the \emph{Farrell-Jones Conjecture} for $RG$ as an input, compare
Example~\ref{exa:FJ}. Modulo the Farrell-Jones Conjecture, we will give a complete answer to
Question~\ref{que:main-question} for instance in the case of Hochschild and  cyclic homology, when the coefficient
ring $R$ is an algebraic number field $F$ or its ring of integers $\calo_F$. We will also give
partial results for  periodic cyclic and negative cyclic homology.

All detection results are obtained by using only the Dennis trace with values in $\HH_{\ast}^{\otimes_{k}} ( RG )$,
whereas all vanishing results hold even for the trace with values
in $\HN_{\ast}^{\otimes_{\IZ}} ( RG )$, which, in view of diagram~\eqref{trace-maps}, can be viewed as the best among
the considered trace invariants.
(Note that for a $k$-algebra $R$ every homomorphism $k^{\prime} \to k$ of commutative rings leads to a homomorphism
$\HN_{\ast}^{\otimes_{k^{\prime}}} ( R ) \to \HN_{\ast}^{\otimes_{k}} ( R )$. Similar for Hochschild, cyclic and periodic
cyclic homology.)
We have no example where the extra effort that goes into the construction of the variants with values in
cyclic, periodic cyclic or negative cyclic homology yields more information about $K_{\ast}(RG)\otimes_{\IZ} \IQ$ than
one can obtain by Hochschild homology, see also Remark~\ref{rem:Chern-characters} and \ref{remark-LRRV} below.

We will now explain our main results. We introduce some notation.

\begin{notation}
Let $G$ be a group and $H$ a subgroup. We write $\langle g \rangle$ for the cyclic
subgroup generated by $g \in G$. We denote by $\conjclass{g}$ and by $(H)$ the conjugacy class of $g$ respectively
of $H$ in $G$. Let $\con G$ be the set of conjugacy classes of elements of $G$.
The set of conjugacy classes of finite cyclic subgroups of $G$ will be denoted by $( \FCyc )$.

Let $Z_G H$ and $N_G H$
denote the centralizer and the normalizer of  $H$ in $G$, respectively.
The \emph{Weyl group}
$W_G H$ is defined as the quotient $N_G H \big/ H \cdot Z_G H$ and coincides for an abelian subgroup $H$ with
$N_G H / Z_G H$.
\end{notation}

Let $C$ be a finite cyclic group. We will define in \eqref{theta_C} an idempotent
$\theta_C \in A(C) \otimes_{\IZ} \IQ$ in the rationalization of the \emph{Burnside ring} $A(C)$ of $C$.
Since there is a natural action of $A(C)$ on $K_{\ast}(RC)$, we obtain a corresponding direct summand
$$\theta_C\big(K_{\ast}(RC) \otimes_{\IZ} \IQ\big) \; \subseteq\; K_{\ast}(RC) \otimes_{\IZ} \IQ\,.$$
In Lemma~\ref{lem:module-over-burnside}, we prove that $\theta_C ( K_{\ast} (RC ) \otimes_{\IZ} \IQ )$ is isomorphic to
the \emph{Artin defect}
\begin{eqnarray*}
\coker \left( \bigoplus_{D \lneqq C} \ind_D^C \colon \bigoplus_{D \lneqq C}
K_{\ast} ( RD ) \otimes_{\IZ} \IQ \to
K_{\ast} ( RC ) \otimes_{\IZ} \IQ
\right),
\end{eqnarray*}
which measures the part of $K_{\ast} ( RC ) \otimes_{\IZ} \IQ$ which is not obtained by induction from proper subgroups
of $C$.

The conjugation action of $N_GC$ on $C$ induces an action of the Weyl group $W_GC = N_GC/Z_GC$
on $K_{\ast} ( RC ) \otimes_{\IZ} \IQ$ and thus on $\theta_C(K_{\ast}(RG) \otimes_{\IZ} \IQ)$.
There is an obvious $W_GC$-action on $BZ_GC = Z_GC\backslash EN_GC$. These actions are understood
in the following statement.

\begin{theorem}[Main Detection Result] \label{the:main-detection-result}
\mbox{}\\
Let $G$ be a group, $k$ a commutative ring and $R$ a $k$-algebra. Suppose that the underlying ring of $R$
is from the following list\,:
\begin{enumerate}
\item a finite dimensional semi\-simple algebra $R$ over a field $F$ of characteristic zero;
\item a commutative complete local domain $R$ of characteristic zero;
\item a commutative Dedekind domain $R$ in which the order of every finite cyclic subgroup of $G$ is invertible and whose quotient field is an algebraic number field.
\end{enumerate}
Then there exists an injective homomorphism
\begin{eqnarray} \label{main-injective-map}
\hspace{3em}
\bigoplus_{(C) \in (\FCyc)}
H_{\ast}(BZ_GC;\IQ) \otimes_{\IQ[W_GC]} \theta_C \big(  K_0(RC) \otimes_{\IZ} \IQ \big)
\to K_{\ast} ( RG ) \otimes_{\IZ } \IQ
\end{eqnarray}
whose image is detected by the Dennis trace map
\begin{eqnarray} \label{Dennis-trace}
\dtr \colon K_{\ast} ( RG ) \otimes_{\IZ} \IQ \to \HH^{\otimes_{k}}_{\ast}( RG ) \otimes_{\IZ } \IQ\,,
\end{eqnarray}
in the sense that the composition of the map \eqref{main-injective-map} with $\dtr$ is injective.
Also the composition with the map to $\HC^{\otimes_{k}}_{\ast}( RG ) \otimes_{\IZ } \IQ$ remains
injective.
\end{theorem}

Examples of rings $R$ appearing in the list of Theorem~\ref{the:main-detection-result} are\,:
\begin{itemize}

\item fields of characteristic zero;

\item the group ring $FH$ of a finite group $H$ over a field $F$ of characteristic zero;

\item the ring $ \IZ\widehat{_p}$ of $p$-adic integers;

\item for the given $G$, the localization $S^{-1}\calo_F$ of the ring of integers $\calo_F$ in an
algebraic number field $F$, for instance $S^{-1}\IZ$, where $S$ is the multiplicative set generated by the
orders of all finite cyclic subgroups of $G$.

\end{itemize}

Depending on the choice of the coefficient ring $R$, the description of the source of the map \ref{main-injective-map}
can be simplified. We mention two examples. Let $\IQ_\infty$ be the
field obtained from $\IQ$ by adjoining all roots of unity.

\begin{theorem}[Detection Result for $\IQ$ and $\IC$ as coefficients]\label{the:positive for IQ and IC}
\mbox{}\\
For every group $G$, there exist injective homomorphisms
\begin{eqnarray*}
\bigoplus_{(C) \in ( \FCyc )} H_{\ast}( B N_G C;\IQ) & \!\!\!\to\!\!\! & K_{\ast} ( \IQ G ) \otimes_{\IZ } \IQ,
\\[.2em]
\bigoplus_{\conjclass{g} \in \con G , |g| < \infty} H_{\ast}( B Z_G \langle g \rangle;\IQ_\infty) & \!\!\!\to\!\!\! &
K_{\ast} ( \IC G ) \otimes_{\IZ } \IQ_\infty\,.
\end{eqnarray*}
The image of these maps is detected by the Dennis trace map with $\IQ$ and $\IC$ as ground ring, respectively.
The coefficient field $\IQ$ (resp.\ $\IC$) can be replaced by any field of characteristic zero (resp.\ any
field containing $\IQ_\infty$).
\end{theorem}

Theorem~\ref{the:positive for IQ and IC} for $\IQ_\infty$ and $\IC$ as
coefficient fields is the main result of the paper by Matthey~\cite{Matthey(2001)}.
The techniques there are based on so-called \emph{delocalization}
and the computation of the Hochschild homology and of the cyclic homology of group
rings with commutative coefficient rings containing $\IQ$ (see
\cite[Section~9.7]{Weibel(1994)} and \cite{Burghelea(1985)}). They
are quite different from the ones used in the present paper and are
exactly suited for the cases studied there and do not seem to be
extendable to the situations considered here. Both maps appearing in
Theorem~\ref{the:positive for IQ and IC} are optimal in the
sense of Theorem~\ref{the:optimal-result} and of
Theorem~\ref{the:Another-optimal-result} below, provided that the
Farrell-Jones Conjecture holds rationally for $K_{\ast}(\IQ G)$ and
$K_{\ast}(\IC G)$ respectively.

The Main Detection Theorem~\ref{the:main-detection-result}  is obtained by studying the following commutative diagram\,:
\begin{eqnarray} \label{diagram-first-appearance}
\vcenter{\xymatrix@C=4.5em{
H_{\ast}^G ( \underline{E}G ; \JK R ) \ar[d]_{H_{\ast}^G( \underline{E}G ; \Jdtr )}  \ar[r]^-{\assembly} & K_{\ast} ( RG ) \ar[d]^{\dtr} \\
\; H_{\ast}^G ( \underline{E}G ; \JH \JH^{\otimes_{\IZ}} R ) \; \ar[r]^-{\assembly} & \HH_{\ast}^{\otimes_{\IZ}} ( RG )
}}
\end{eqnarray}
Here, the horizontal arrows are \emph{generalized assembly maps} for $K$-theory and Hochschild homology respectively, and
the left vertical arrow is a suitable version of the Dennis trace map.
The $G$-space $\underline{E}G$ is a model for the so-called \emph{classifying space for proper $G$-actions}.
Moreover, $H_{\ast}^G( - ; \JK R)$ and $H_{\ast}^G( - ;  \JH \JH^{\otimes_{\IZ}} R )$ are certain
$G$-homology theories. We will explain the diagram in more detail in
Section~\ref{sec:Outline-of-the-method}. We will prove that the lower horizontal arrow in \eqref{diagram-first-appearance}
is split injective, see Theorem~\ref{the:splitHochschild}.
In fact, Theorem~\ref{the:splitHochschild} gives a complete picture of the generalized assembly map for
Hochschild and cyclic homology.
We will also compute the left-hand vertical arrow after rationalization, compare Theorem~\ref{the:Chern} and
Propositions \ref{prop:compute-theta-null}, \ref{prop:comp-Z} and \ref{prop:compute-higher-HH-HC}. According
to this computation, the left-hand side in \eqref{main-injective-map} is a direct summand in
$H_{\ast}( \underline{E}G ; \JK R ) \otimes_{\IZ} \IQ$ on which, for $R$ as in Theorem~\ref{the:main-detection-result},
the map
\begin{eqnarray} \label{dtr-homology}
H^G_{\ast} ( \underline{E}G ; \JK R ) \otimes_{\IZ} \IQ \to H^G_{\ast} ( \underline{E}G ;
\JH \JH^{\otimes_{\IZ}} R ) \otimes_{\IZ} \IQ
\end{eqnarray}
is injective. This will prove Theorem~\ref{the:main-detection-result}. Now, suppose that $R$ is as in case $(i)$ of
Theorem~\ref{the:main-detection-result}, with $F$ a number field. Then, it turns out that the map \eqref{dtr-homology}
vanishes on a complementary summand. According to the Farrell-Jones Conjecture for $K_{\ast}(RG)$, the upper
horizontal arrow in \eqref{diagram-first-appearance} should be an isomorphism (this uses that  $R$ is a regular
ring with $\IQ \subseteq R$). Combining these facts, we will deduce the following result.

\begin{theorem}[Vanishing Result for Hochschild and cyclic homology] \label{the:optimal-result}
\mbox{}\\
Let $G$ be a group, $F$ an algebraic number field, and $R$ be finite dimensional semi\-simple $F$-algebra.
Suppose that for some $n\geq 0$, the Farrell-Jones Conjecture holds rationally for $K_{n}(RG)$, see
Example~\ref{exa:FJ} below.

Then Theorem~\ref{the:main-detection-result} is optimal for the Hochschild homology trace invariant, in the sense that
the Dennis trace map
\begin{eqnarray} \label{dennis-with-z-coefficients}
\dtr\colon K_{n} ( RG ) \otimes_{\IZ} \IQ \to \HH^{\otimes_{\IZ}}_{n} ( RG )\otimes_{\IZ} \IQ
\end{eqnarray}
vanishes on a direct summand that is complementary to the image of the injective map \eqref{main-injective-map}
in degree $n$. Consequently, also the trace taking values in rationalized cyclic homology
$\HC^{\otimes_{\IZ}}_{n} ( RG ) \otimes_{\IZ} \IQ$ vanishes on this complementary summand.
\end{theorem}

One might still hope that the refinements of the Dennis trace map with values in periodic cyclic or negative
cyclic homology detect more of the rationalized algebraic $K$-theory of $RG$.
But one can show that this is not the case if one additionally assumes
a finiteness condition on the classifying space $\underline{E}G$.
Recall that the $G$-space $\underline{E}G$ is called \emph{cocompact} if the orbit space
$G\backslash \underline{E}G$ is compact, in other words, if it consists of finitely many $G$-equivariant cells.
Cocompact models for $\underline{E}G$ exist for many interesting groups
$G$ such as discrete cocompact subgroups of virtually connected Lie groups, word-hyperbolic groups, arithmetic subgroups of
a semi\-simple connected $\IQ$-algebraic group, and
mapping class groups (see for instance \cite{Lueck(2004h)}).

\begin{theorem}[Vanishing Result for periodic and negative cyclic homology] \label{the:Another-optimal-result}
\mbox{}\\
Let $F$ be an algebraic number field, and $R$ a finite dimensional semi\-simple $F$-algebra.
Suppose that for some $n\geq 0$, the Farrell-Jones Conjecture holds rationally for $K_{n} ( RG )$. Suppose
further that there exists a cocompact model for the classifying space for proper $G$-actions \underline{E}G.

Then also the refinements of the Dennis trace with values in $\HP^{\otimes_{\IZ}}_{n} ( RG ) \otimes_{\IZ} \IQ$ and
in $\HN^{\otimes_{\IZ}}_{n} ( RG ) \otimes_{\IZ} \IQ$ vanish on a direct summand which is complementary to the image
of the injective map \eqref{main-injective-map} in degree $n$.
\end{theorem}

The next result is well-known. It shows in particular that the rational group homology
$H_{\ast}(BG;\IQ)$ is contained in $K_{\ast}(RG) \otimes_{\IZ} \IQ$ for all commutative rings $R$ of characteristic
zero.

\begin{theorem}[Detection Result for commutative rings of characteristic zero]
\label{the:detection-result-for-comm.-rings-of-char.-zero}
\mbox{}\\
Let $R$ be a ring such that the canonical ring homomorphism $\IZ \to R$ induces an
injection $\HH_0^{\otimes_{\IZ}}(\IZ) = \IZ \into \HH_0^{\otimes_{\IZ}}(R) = R\big/[R,R] $, for instance a
commutative ring of characteristic zero.

Then, for any group $G$, there exists an injective homomorphism
\begin{eqnarray} \label{injective-for-z-coefficients}
H_{\ast} ( BG ; \IQ ) \to K_{\ast} (RG) \otimes_{\IZ } \IQ
\end{eqnarray}
whose composition with the Dennis trace map \eqref{Dennis-trace}
is injective for every choice of a ground ring $k$ such that $R$ is a $k$-algebra.
The corresponding statement holds with Hochschild homology replaced by cyclic homology.
\end{theorem}

Special cases of this result are treated for example in \cite[Proposition~6.3.24 on page~366]{Rosenberg(1994)}.

According to the Farrell-Jones Conjecture, the image of \eqref{injective-for-z-coefficients} should only be a
very small part of the rationalized $K$-theory of $RG$.
The following result illustrates that, for certain coefficient rings, including $\IZ$, one cannot expect to detect
more by linear traces than achieved in Theorem~\ref{the:detection-result-for-comm.-rings-of-char.-zero}.

\begin{theorem}[Vanishing Result for integral coefficients]  \label{the:optimal-Z-coeff}
\mbox{}\\
Let $S^{-1}\calo_F$ be a localization of a ring of integers $\calo_F$  in an algebraic number field $F$
with respect to a (possibly empty) multiplicatively closed subset $S$.
Assume that no prime divisor of the order $|H|$ of a nontrivial finite subgroup $H$ of $G$ is invertible in $S^{-1} \calo_F$.
Suppose that for some $n\geq 0$, the Farrell-Jones Conjecture holds rationally for
$K_n(S^{-1} \calo_F[G])$.

Then the Dennis trace \eqref{dennis-with-z-coefficients} vanishes on a summand in
$K_n (S^{-1} \calo_F[G])\otimes_{\IZ} \IQ$
which is complementary to the image of
the map \eqref{injective-for-z-coefficients} in degree $n$.
Consequently, the analogous statement holds for the trace with values in
$\HC^{\otimes_{\IZ}}_{n} (S^{-1} \calo_F[G]) \otimes_{\IZ} \IQ$.
\end{theorem}

The most interesting case in Theorem~\ref{the:optimal-Z-coeff} is $R = \IZ$.
We remark that rationally, the Farrell-Jones Conjecture for $K_{\ast}( \IZ G )$
is known in many cases, for example for every subgroup $G$ of a discrete cocompact subgroup of a
virtually connected Lie group \cite{Farrell-Jones(1993a)}. For a survey of known results about the
Farrell-Jones Conjecture, we refer the reader to \cite{Lueck-Reich(2003b)}.

\begin{remark} \label{rem:Chern-characters}
There are further trace invariants (or Chern characters) given by maps $\ch_{n,r}\colon K_{n}
( RG ) \to \HC^{\otimes_{k}}_{n+2r} ( RG )$, for fixed $n$, $r \geq 0$, see~\cite[8.4.6 on page 272 and 11.4.3 on page 371]{Loday(1992)}. This
will however produce no new detection results in the spirit of the above statements, since there is a commutative
diagram
\[
\xymatrix@R=3em@C=.5em{
& & & & K_{n} ( RG ) \ar[dllll]_{\ntr} \ar[dl] \ar[dr]_(.4){\ch_{n,r}\!\!\!\!\!} \ar[drrrr] & & & & \\
\HN_{n}^{\otimes_k} ( RG ) \ar[rrr] & & & \HP_{n}^{\otimes_k} ( RG ) \ar[rr]^-{\overline{S}} & & \HC_{n+2r}^{\otimes_k} ( RG )
\ar[rrr]^-{S^{r}} & & & \HC_{n}^{\otimes_k} ( RG ).
         }
\]
\end{remark}

\begin{remark} \label{remark-LRRV}
In \cite{Boekstedt-Hsiang-Madsen(1993)}, B{\"o}kstedt, Hsiang and Madsen
define the \emph{cyclotomic trace}, a map out of $K$-theory which takes values in
\emph{topological cyclic homology}. The cyclotomic trace map can be thought of as an even more elaborate refinement of the
Dennis trace map. In contrast to the Dennis trace, it seems that the cyclotomic trace has the potential to detect
almost all of the rationalized $K$-theory of an integral group ring. This question is investigated in detail in
\cite{Lueck-Reich-Rognes-Varisco(2005)}.
\end{remark}

%\begin{remark}  \label{remark-about-R-not-free-over-k}
%%By $\HH^{\otimes_k}( R )$ we denote the Hochschild homology defined via the standard
%Hochschild complex using tensor products over $k$ (compare \cite[Chapter~1]{Loday(1992)}).
%We do not assume that $R$ (or $RG$) is flat
%over $k$. Note that if $R$ is not flat over $k$ then depending on the point of view this may be a
%too ``naive'' definition of Hochschild homology, because it does not have the more conceptional interpretation
%as $\Tor_k^{R\otimes_k R^{\op}}( R ,R )$, compare \cite[Proposition~1.1.13 on page~12]{Loday(1992)}.
%This defect can however be repaired using the framework of relative homological algebra, compare
%\cite[Lemma~9.1.3 on page~302]{Weibel(1994)}.
%\end{remark}

\medskip

The paper is organized as follows\,:

\smallskip

\begin{tabular}{ll}
\ref{sec:Outline-of-the-method}. & Outline of the method
\\
\ref{sec:Proofs}. & Proofs
\\
\ref{sec:trace-for-finite-cyclic}. & The trace maps for finite cyclic groups
\\
\ref{sec:Notation}. & Notation and general machinery
\\
\ref{sec:trace-maps}. & The trace maps
\\
\ref{sec:mackey}. & Equivariant homology theories, induction and Mackey structures
\\
\ref{sec:eval-chern}. & Evaluating the equivariant Chern character
\\
\ref{sec:comparing-different-models}. & Comparing different models
\\
\ref{sec:Splitting-assembly-maps}. & Splitting assembly maps
\\
 & References
\end{tabular}

%%%%%%%%%%%%%%%%%%%%%%%%%%%%%%%%%%%%%%%%%%%%%%%%%%%%%%%%%%%%%%%%%%%%%%%%%%%%%%%%%%%%%%%%%%%%%%%%%%%%%%
%%%%%%%%%%%%%%%%%%%%%%%%%%%%%%%%%%%%%%%%%%%%%%%%%%%%%%%%%%%%%%%%%%%
\section{Outline of the method}\label{sec:Outline-of-the-method}

This paper is concerned with comparing generalized assembly maps for $K$-theory, via the Dennis trace or its refinements,
with generalized assembly maps for Hochschild homology, for cyclic, periodic cyclic or negative cyclic homology.
Before we explain the general strategy behind our results we briefly explain the concept of a generalized assembly map; for more details the reader is referred
to \cite{Davis-Lueck(1998)} and \cite[Section~2 and 6]{Lueck-Reich(2003b)}.

A \emph{family of subgroups} of a given group $G$ is a non-empty collection of subgroups which is closed under
conjugation and finite intersections. Given a family $\calf$ of subgroups, there always exists a $G$-$CW$-complex
$\EGF{G}{\calf}$ all of whose isotropy groups lie in $\calf$ and which has the property that for all $H \in \calf$, the
fixed subspace $\EGF{G}{\calf}^H$ is a contractible space. A $G$-$CW$-complex with these properties is unique up to
$G$-homotopy because it receives a $G$-map from every $G$-$CW$-complex all whose isotropy groups lie in $\calf$ and this $G$-map
is unique up to $G$-homotopy. If $\calf= \Fin$ is the \emph{family of finite subgroups}, then one often writes $\underline{E}G$
for $\EGF{G}{\Fin}$. For a survey on these spaces, see for instance \cite{Lueck(2004h)}.

%We also write $(\calf)$ for the
%\emph{set of conjugacy classes $(H)$ of subgroups $H\in\calf$}.

Let $\Or G$ denote the \emph{orbit category} of $G$. Objects are the homogenous spaces $G/H$ considered as left $G$-spaces, and
morphisms are $G$-maps. A functor $\JE$, from the orbit category $\Or G$ to the category of spectra, is called an
\emph{$\Or G$-spectrum}. Each $\Or G$-spectrum $\JE$ gives rise to a
\emph{$G$-homology theory} $H^G_{\ast}( - ; \JE )$, compare~\cite[Section~6]{Lueck-Reich(2003b)} and the beginning of
Section~\ref{sec:mackey} below. Given $\JE$ and a family $\calf$ of subgroups of $G$, the so-called  \emph{generalized assembly map}
\begin{eqnarray} \label{assembly-map}
\xymatrix@C=4em{
H^G_{\ast} \big(\EGF{G}{\calf}; \JE \big) \ar[r]^-{\assembly} & H^G_{\ast}( \pt ; \JE )
               }
\end{eqnarray}
is merely the homomorphism induced by the map $\EGF{G}{\calf} \to \pt$. The group
$H_{\ast}^G( \pt ; \JE)$ can be canonically identified with $\pi_{\ast}( \JE( G/G ))$.

\begin{example}[The Farrell-Jones Conjecture] \label{exa:FJ}
\mbox{}\\
Given an arbitrary ring $R$ and an arbitrary group $G$, there exists a \emph{non-connective $K$-theory $\Or G$-spectrum},
denoted by $\JK^{\minusinfinity}R(?)$, such that there is a natural isomorphism
$$
\pi_n \big(\JK^{\minusinfinity}R ( G/H ) \big) \cong K_n ( RH )
$$
for all $H\leq G$ and all $n \in \IZ$, compare \cite[Theorem~6.9]{Lueck-Reich(2003b)}.
The \emph{Farrell-Jones Conjecture} for $K_{n}(RG)$, \cite[1.6 on page~257]{Farrell-Jones(1993a)},
predicts that the generalized assembly map
\[
\xymatrix@C=4em{
H_n^G \big(\EGF{G}{\VCyc};\JK^{\minusinfinity} R \big) \ar[r]^-{\assembly} &
H_n^G(\pt;\JK^{\minusinfinity} R) \cong K_n(RG)
                 }
\]
is an isomorphism. Here $\VCyc$ stands for the \emph{family of all virtually cyclic subgroups} of~$G$.
A group is called virtually cyclic if it contains a cyclic subgroup of finite index.

%We also point out that there is a map of $\Or G$-spectra $\JK R(?) \to \JK^{\minusinfinity} R(?)$
%(the $(-1)$-connective covering map), which is natural in $R$, and which induces an isomorphism
%\[
%\pi_{n}\big(\JK R ( G/H ) \big) \xrightarrow{\cong} \pi_{n}\big(\JK^{\minusinfinity} R ( G/H ) \big)\,,
%\]
%that corresponds to the identity of $K_n(RH)$, for every $H\leq G$ and every $n\geq 0$.
\end{example}

In Section~\ref{sec:trace-maps}, we will construct the following commutative diagram of connective
$\Or G$-spectra and maps (alias natural transformations) between them\,:
\begin{eqnarray} \label{spec-level-trace-maps}
\vcenter{\xymatrix@C=3.2em{
& \JH \JN^{\otimes_k} R \ar[r] \ar[d]^{\Jh} & \JH \JP^{\otimes_k} R \ar[d]  \\
\JK R \ar[r]^-{\Jdtr} \ar[ur]^-{\Jntr} & \JH \JH^{\otimes_k} R \ar[r] & \JH \JC^{\otimes_k} R.
        }}
\end{eqnarray}
Decisive properties of these constructions are that for all $n \geq 0$, we have natural isomorphisms
\begin{eqnarray} \label{relate-to-classical}
\renewcommand{\arraystretch}{1.4}
\begin{array}{lcl}
\pi_n\big(\JK R ( G/H )\big) & \cong  & K_n ( RH )
\\
\pi_n\big( \JH \JH^{\otimes_k } R ( G/H ) \big) & \cong & \HH_n^{\otimes_k} ( RH )
\\
 \pi_n\big( \JH \JC^{\otimes_k } R ( G/H ) \big) & \cong & \HC_n^{\otimes_k} ( RH )
\\
\pi_n\big( \JH \JP^{\otimes_k } R ( G/H ) \big)  & \cong  &
\HP_n^{\otimes_k} ( RH )
\\
\pi_n\big( \JH \JN^{\otimes_k } R ( G/H ) \big)  & \cong  &
\HN_n^{\otimes_k} ( RH )
\end{array}
\end{eqnarray}
and all negative homotopy groups vanish. Note that we need to distinguish between the non-connective version
$\JK^{-\infty}R$ and the connective version $\JK R$.
Under the identifications above, the maps of $\Or G$-spectra in \eqref{spec-level-trace-maps}
evaluated at an orbit $G/H$ induce, on the level of homotopy groups, the maps in \eqref{trace-maps}
with $R$ replaced by the corresponding group ring $RH$.

\begin{remark} We found it technically convenient to work, at the level of spectra, with the connective versions
of periodic cyclic and negative cyclic homology. Since we are mainly interested in the trace maps (whose source
will be the connective $K$-theory spectrum), we do not lose any information.
\end{remark}

Since the assembly map \eqref{assembly-map} is natural in the functor $\JE$, we obtain, for each family of
subgroups $\calf$ of a group $G$ and for each $n\geq 0$, the commutative diagram
\begin{eqnarray} \label{diagram-assembly-dtr}
\mbox{\phantom{xxxxxx}}\vcenter{\xymatrix@C=4.5em{
H_n^G \big( \EGF{G}{\calf};\JK R \big) \ar[r]^-{\assembly} \ar[d]_{H_n^G( \EGF{G}{\calf} ; \Jntr)} &
H_n^G (\pt;\JK R) \cong K_n(RG) \ar@<-2.5em>[d]^{\ntr} \\
H_n^G \big( \EGF{G}{\calf}; \JH \JN ^{\otimes_k} R \big) \ar[r]^-{\assembly} \ar[d]_{H_n^G( \EGF{G}{\calf} ; \Jh)} &
H_n^G( \pt ; \JH \JN^{\otimes_k} R ) \cong \HN_n^{\otimes_k} ( RG ) \ar@<-2.5em>[d]^{\h} \\
H_n^G \big( \EGF{G}{\calf};\JH \JH^{\otimes_k}R \big) \ar[r]^-{\assembly} &
H_n^G(\pt;\JH \JH^{\otimes_k}R) \cong \HH_n^{\otimes_k}(RG).
         }}
\end{eqnarray}
The vertical compositions are the corresponding versions of the Dennis trace map.

\smallskip

Our investigation relies on two main ingredients. The first ingredient are splitting and isomorphism
results for the assembly maps of Hochschild and cyclic type.

\begin{theorem} [The Isomorphism Conjecture for $\HH$ and $\HC$] \label{the:splitHochschild}
\mbox{}\\
Let $k$ be a commutative ring, $R$ a $k$-algebra, and $G$ a group. Then the generalized Hochschild homology
assembly map
\[
\xymatrix@C=5em{
\;H_{\ast} \big( \EGF{G}{\calf}; \JH \JH^{\otimes_k} R \big)\; \ar[r]^-{\assembly} &
H_{\ast} ( \pt ; \JH \JH^{\otimes_k} R )\cong \HH_{\ast}^{\otimes_k}(RG)
               }
\]
is split injective for every family $\calf$.
If $\calf$ contains the family of all (finite and infinite) cyclic subgroups, then
the map is an isomorphism. The analogous statement holds for $\JH \JC$ in
place of $\JH \JH$.
\end{theorem}

The fact that the definition of periodic cyclic and of negative cyclic homology involves certain inverse limit processes
prevents us from proving the analogous result in
these cases without assumptions on the group $G$. But we still have the following statement.

\begin{addendum} [Splitting Results for the $\HP$ and $\HN$-assembly maps] \label{add:split-HN}
\mbox{}\\
Suppose that there exists a cocompact  model for the classifying space $\EGF{G}{\calf}$. Then the statement of
Theorem~\ref{the:splitHochschild} also holds for $\JH \JP$ and $\JH \JN$ in place of $\JH \JH$.
\end{addendum}

The proofs of Theorem~\ref{the:splitHochschild} and Addendum~\ref{add:split-HN} are presented in
Section~\ref{sec:Splitting-assembly-maps}.

\begin{remark}
We do not know any non-trivial example where the isomorphism statement in Addendum~\ref{add:split-HN}
applies, i.e.\ where $\calf$ contains all (finite and infinite) cyclic groups and where, at the same
time, $\EGF{G}{\calf}$ has a cocompact model.
\end{remark}

The second main ingredient of our investigation is the rational computation
of equivariant homology theories from \cite{Lueck(2002b)}.
For varying $G$, our $G$-homology theories
like $\calh^G_\ast(-)=H_\ast^G( - ; \JK R )$ or $\calh^G_\ast(-)=H_\ast^G(- ; \JH \JH^{\otimes_k} R)$
are linked by a so-called \emph{induction structure} and form an \emph{equivariant homology theory} in the sense of
\cite{Lueck(2002b)}. Moreover, these homology theories admit a \emph{Mackey structure}.
In Section~\ref{sec:mackey}, we review these notions and explain some general principles which allow us to
verify that $G$-homology theories like the ones we are interested in indeed admit induction and Mackey structures.
In particular, Theorems~0.1 and 0.2 in \cite{Lueck(2002b)} apply and yield an explicit computation of
$\calh^G_\ast ( \underline{E}G ) \otimes_{\IZ} \IQ$. In Section~\ref{sec:eval-chern}, we review this
computation and discuss a simplification which occurs in the case of $K$-theory, Hochschild,
cyclic, periodic cyclic and negative cyclic homology, due to the fact that in all these special cases,
we have additionally a module structure over the Swan ring.

In order to state the result of this computation, we introduce some more notation. For a finite group $G$, we
denote by $A(G)$ the \emph{Burnside ring} which is additively generated by isomorphism classes of finite transitive
$G$-sets.
Let $( \sub G )$ denote the set of conjugacy classes of subgroups of $G$.

The counting fixpoints ring homomorphism
\begin{eqnarray} \label{def-chi-G}
\chi_G \colon A(G) \to \prod_{(\sub G)} \IZ\,
\end{eqnarray}
which is induced by sending  a $G$-set $S$ to $( |S^H| )_{(H)}$ becomes
an isomorphism after rationalization,
compare \cite[page~19]{Dieck(1987)}. For a finite cyclic group $C$, we consider the idempotent
\begin{eqnarray} \label{theta_C}
\theta_C \; = \; (\chi_C \otimes_{\IZ} \IQ)^{-1} \big((\delta_{CD})_{D}\big) & \in & A(C) \otimes_{\IZ} \IQ\,,
\end{eqnarray}
where $(\delta_{CD})_D\in \prod_{\sub C} \IQ$ is given by $\delta_{CC}=1$ and $\delta_{CD}= 0$ if $D \neq C$.

\smallskip

Recall that $K_{\ast} ( RC )$ and similarly Hochschild, cyclic, periodic cyclic and negative cyclic homology
of $RC$ are modules over the Burnside ring $A(C)$. The action of a $C$-set $S$ is in all cases induced from taking
the tensor product over $\IZ$ with the corresponding permutation module $\IZ S$.
In Lemma~\ref{lem:module-over-burnside} below, we prove that $\theta_C ( K_{\ast} ( RC ) \otimes_{\IZ} \IQ )$
is isomorphic to the $\IQ$-vector space
\begin{eqnarray} \label{eq:SC-ist-Artin-defekt}
\hspace*{3em}
\coker \left( \bigoplus_{D \lneqq C} \ind_D^C \colon \bigoplus_{D \lneqq C}
K_{\ast} ( RD ) \otimes_{\IZ} \IQ \to
K_{\ast} ( RC ) \otimes_{\IZ} \IQ
\right),
\end{eqnarray}
which is known as the \emph{Artin defect} of $K_{\ast} ( RC ) \otimes_{\IZ} \IQ$.

In Section~\ref{sec:eval-chern} we establish the following result.

\begin{theorem} \label{the:Chern}
For each $n\geq 0$, the following
diagram commutes and the arrows labelled $\ch^G$ are isomorphisms\,:
\begin{eqnarray*}
\xymatrix@C=-10mm{
\smash[b]{{\displaystyle {\bigoplus_{\substack{p,q \geq 0\\p+q = n}}}\; \displaystyle
\bigoplus_{(C) \in (\FCyc)}}}
H_p(BZ_GC;\IQ) \otimes_{\IQ[W_GC]} \theta_C \big(  K_q (RC) \otimes_{\IZ} \IQ \big)
\ar[dr]^-{\ch^G}_-{\cong}
\ar[dd]_{\dtr_\ast}%_-{\id \otimes_{\IQ[W_GC]} \theta_C (\id \otimes_{\IZ} \dtr_q)}
& \\
& H_n^G(\underline{E}G;\JK R ) \otimes_{\IZ} \IQ
 \ar[dd]^{H_n^G( \underline{E}G ; \Jdtr) \otimes_{\IZ} \IQ} %_-{\id \otimes_{\IZ} H_n^G\underline{E}G; \dtr)}
\\
\smash[b]{{\displaystyle
{\bigoplus_{\substack{p,q \geq 0\\p+q = n}}}\; \displaystyle \bigoplus_{(C) \in (\FCyc )}}}
H_p(BZ_GC;\IQ) \otimes_{\IQ[W_GC]} \theta_C \big(  \HH_q^{\otimes_k}(RC) \otimes_{\IZ} \IQ \big)
\ar[dr]^-{\ch^G}_-{\cong}
& \\
& H_n^G(\underline{E}G;\JH \JH^{\otimes_k} R) \otimes_{\IZ} \IQ.
}
\end{eqnarray*}
The left-hand vertical arrow is induced by the Dennis trace maps for finite cyclic groups and
respects the double direct sum decompositions. The right-hand vertical arrow is induced by
the $\Or G$-spectrum Dennis trace $\Jdtr$, compare \eqref{spec-level-trace-maps}.
There are similar diagrams and isomorphisms
corresponding to each of the other maps in diagram \eqref{spec-level-trace-maps}.
\end{theorem}

\begin{remark} \label{remark-nonconn-conn}
The $(-1)$-connected covering map of $\Or G$-spectra $\JK R \to \JK^{-\infty} R$
induces for every orbit $G/H$ an isomorphism
\[
\pi_n ( \JK R ( G/H) ) \to \pi_n ( \JK^{-\infty} R ( G/H) )
\]
if $n \geq 0$. The source is trivial for $n < 0$. This map induces
the following commutative diagram.
\begin{eqnarray*}
\xymatrix@C=-10mm{
\smash[b]{{\displaystyle {\bigoplus_{\substack{p,q \geq 0\\p+q = n}}}\; \displaystyle
\bigoplus_{(C) \in (\FCyc)}}}
H_p(BZ_GC;\IQ) \otimes_{\IQ[W_GC]} \theta_C \big(  K_q (RC) \otimes_{\IZ} \IQ \big)
\ar[dr]^-{\ch^G}_-{\cong}
\ar[dd]%_-{\id \otimes_{\IQ[W_GC]} \theta_C (\id \otimes_{\IZ} \dtr_q)}
& \\
& H_n^G(\underline{E}G;\JK R ) \otimes_{\IZ} \IQ
 \ar[dd]%_-{\id \otimes_{\IZ} H_n^G\underline{E}G; \dtr)}
\\
\smash[b]{{\displaystyle
{\bigoplus_{\substack{p,q \in \IZ \\p+q = n}}}\; \displaystyle \bigoplus_{(C) \in (\FCyc )}}}
H_p(BZ_GC;\IQ) \otimes_{\IQ[W_GC]} \theta_C \big(  K_q (RC) \otimes_{\IZ} \IQ \big)
\ar[dr]^-{\ch^G}_-{\cong}
& \\
& H_n^G(\underline{E}G;\JK^{\minusinfinity} R) \otimes_{\IZ} \IQ.
}
\end{eqnarray*}
Here the arrows
labelled $\ch^G$ are isomorphisms.
Note the restriction $p$, $q \geq 0$ for the sum in the upper left hand corner.
\end{remark}

\subsection{} \label{subsec:general-strategy}
{\bf General strategy.}
We now explain the strategy behind all the results that appeared in the introduction.
If we combine the diagram appearing in Theorem~\ref{the:Chern} with diagram \eqref{diagram-assembly-dtr},
for each $n\geq 0$, we get a commutative diagram
\begin{eqnarray*} \label{diagram3}
\xymatrix@C=-1mm{
\smash[b]{{\displaystyle {\bigoplus_{\substack{p,q \geq 0\\p+q = n}}}\; \displaystyle
\bigoplus_{(C) \in (\FCyc)}}}
H_p(BZ_GC;\IQ) \otimes_{\IQ[W_GC]} \theta_C \big(  K_q (RC) \otimes_{\IZ} \IQ \big)
\ar[dr]%^{\assembly\circ\ch^{G}}
\ar[dd]_{\dtr_\ast}%_-{\id \otimes_{\IQ[W_GC]} \theta_C (\id \otimes_{\IZ} \dtr_q)}
& \\
& K_n(RG) \otimes_{\IZ} \IQ
 \ar[dd]^{\dtr}
\\
\smash[b]{{\displaystyle
{\bigoplus_{\substack{p,q \geq 0\\p+q = n}}}\; \displaystyle \bigoplus_{(C) \in (\FCyc )}}}
H_p(BZ_GC;\IQ) \otimes_{\IQ[W_GC]} \theta_C \big(  \HH_q^{\otimes_\IZ}(RC) \otimes_{\IZ} \IQ \big)
\ar[dr]
& \\
& \;\HH_n^{\otimes_{\IZ}}(RG) \otimes_{\IZ} \IQ.\;
}
\end{eqnarray*}
Because of Theorem~\ref{the:splitHochschild} and the isomorphism statement in Theorem~\ref{the:Chern} the lower horizontal map
is injective.
There is an analogue of the commutative diagram above, where
the upper row is the same and $\HH$ is replaced by $\HC$ in the bottom row. Also in this case we know that the lower horizontal map
is injective because of
Theorem~\ref{the:splitHochschild} and \ref{the:Chern}.

\smallskip

Observe that $W_G C$ is always a finite group, hence $\IQ[W_G C]$ is a semi\-simple ring, so that every module
over $\IQ W_G C$ is flat and the functor $H_p( B Z_G C ; \IQ) \otimes_{\IQ W_G C} ( - )$
preserves injectivity.

For $q\geq 0$ given, we see that suitable injectivity results
about the maps
\begin{eqnarray} \label{eq:the-map-theta-C}
\theta_C \big( K_q( RC ) \otimes_{\IZ} \IQ \big)  \to \theta_C \big( \HH_q^{\otimes_k} ( RC ) \otimes_{\IZ} \IQ \big)
\end{eqnarray}
for the finite cyclic subgroups $C\leq G$ lead to the proof
of detection results in degree $n$. These maps \eqref{eq:the-map-theta-C} will be studied in
Section~\ref{sec:trace-for-finite-cyclic}.

If $R$ is a regular ring containing $\IQ$, then the family $\VCyc$ of virtually cyclic subgroups
can be replaced by the family $\Fin$ of finite subgroups and the non-connective $K$-theory $\Or G$-spectrum
$\JK^{\minusinfinity} R(?)$ by its connective version $\JK R(?)$ in the statement of the Farrell-Jones Conjecture,
i.e., in this case, the Farrell-Jones Conjecture for $K_n(RG)$, for some $n\in\IZ$, is equivalent to the statement
that the assembly map
$$
\xymatrix@C=4em{
H_n^G(\underline{E}G; \JK R) \ar[r]^-{\assembly} & K_n(RG)
               }
$$
is an isomorphism if $n \geq 0$ and to the statement that $K_n(RG) = 0$ if $n \leq -1$
(see \cite[Proposition~2.14]{Lueck-Reich(2003b)}). As a consequence, the upper horizontal
arrow in the diagram above (where $n\geq 0$) is bijective if the Farrell-Jones Conjecture
is true rationally for $K_n(RG)$.

So for $q\geq 0$ given, we see that suitable vanishing results
about the maps \eqref{eq:the-map-theta-C}
(and about their analogues involving cyclic homology)
combined with the assumption that the Farrell-Jones conjecture holds rationally for $K_n(RG)$ lead to the proof of
vanishing results in degree $n$.

\section{Proofs} \label{sec:Proofs}

Based on the strategy explained in the previous paragraphs we now give the proofs
of the theorems stated in the introduction,
modulo the following results\,: Theorem~\ref{the:splitHochschild} and Addendum~\ref{add:split-HN}
(proved in Section~\ref{sec:Splitting-assembly-maps}); Theorems \ref{the:Chern} (proved in Section~\ref{sec:eval-chern}, using Sections
\ref{sec:Notation}--\ref{sec:mackey});
and the results of
Section~\ref{sec:trace-for-finite-cyclic} (which is self-contained, except for Lemma~\ref{lem:module-over-burnside}
whose proof is independent of the rest of the paper).

\subsection{}
{\bf Proof of Theorem~\ref{the:main-detection-result}.}
After the general strategy~\ref{subsec:general-strategy}, the necessary injectivity result to complete the proof
appears in Proposition~\ref{prop:compute-theta-null} below.
\qed

\subsection{}
{\bf Proof of Theorem~\ref{the:optimal-result}.}
The result follows directly from the general strategy~\ref{subsec:general-strategy} and the vanishing result
stated as Proposition~\ref{prop:compute-higher-HH-HC} below.
\qed

\subsection{}
{\bf Proof of Theorem~\ref{the:Another-optimal-result}.}
The proof is completely analogous to that of Theorem~\ref{the:optimal-result}. The extra condition that there
is a cocompact model for $\underline{E}G$ is only needed to apply Addendum~\ref{add:split-HN} in place of
Theorem~\ref{the:splitHochschild}.
\qed

\subsection{}
{\bf Proof of Theorem~\ref{the:positive for IQ and IC}.}
The next lemma explains why Theorem~\ref{the:positive for IQ and IC} for $\IQ$ as coefficients
follows from Theorem~\ref{the:main-detection-result}. The case of $\IC$ as coefficients is proven
similarly, compare \cite[Example~8.11]{Lueck(2002b)}.
\qed

\begin{lemma} \label{lem:compute-theta}
\begin{enumerate}
\item \label{lem:compute-theta-i} Let $C$ be a finite cyclic group. Then one has
\[
\theta_C \big( K_0 ( \IQ C ) \otimes_{\IZ} \IQ \big) \; \cong \; \IQ
\]
and every group automorphism of $C$ induces the identity on $\IQ$.

\item  \label{lem:compute-theta-ii} For any group $G$ and finite cyclic subgroup $C \leq G$, the map
$$H_{\ast}(BZ_GC;\IQ) \otimes_{\IQ[W_GC]} \IQ \xrightarrow{\cong} H_{\ast}(BN_GC;\IQ)$$
induced by the inclusion $Z_GC \into N_GC$ is an isomorphism. Here $\IQ$ carries the trivial $W_GC$-action.
\end{enumerate}
\end{lemma}
\begin{proof}
\ref{lem:compute-theta-i}
There is a commutative diagram
\[
\xymatrix{
A(C) \otimes_{\IZ} \IQ \ar[d]_{\chi_C \otimes_{\IZ} \IQ}^-{\cong} \ar[r]^-{\cong} & K_0 ( \IQ C )
\otimes_{\IZ} \IQ \ar[d]^-{\cong} \\
\prod_{D \in \sub C} \IQ \ar[r]^-{\cong} & \map( \sub C , \IQ )
         }
\]
Here, the upper horizontal map sends a $C$-set to the corresponding permutation module.
The product in the lower left corner is taken
over the set $\sub C$ of all subgroups of $C$ and the left-hand vertical arrow is given by
sending the class of a $C$-set $S$ to $(|S^D|)_D$ and is an isomorphism, as already mentioned after \eqref{def-chi-G}.
The right-hand vertical map is given by sending a rational representation $V$ to its character, i.e.\ if $d$ generates the subgroup $\langle d \rangle$, then
$\langle d \rangle \mapsto \tr_{\IQ} (d \colon V \to V)$. This map is also an isomorphism, compare
\cite[II.\,\S\,12]{Serre(1977)}.
The lower horizontal map is the isomorphism given by sending $(x_D)_{D \in \sub C}$ to $(D \mapsto x_D)$.
The diagram is natural with respect to automorphisms of
$C$. By definition, $\theta_C \in A(C) \otimes_{\IZ} \IQ$  corresponds to the idempotent
$(\delta_{CD})_D$ in the lower left-hand corner. Now, the result follows.

\smallskip

\ref{lem:compute-theta-ii} This follows from the Lyndon-Hochschild-Serre spectral
sequence of the fibration $BZ_GC \to BN_GC \to BW_GC$ and from the fact that,
the group $W_GC$ being finite, for any $\IQ[W_GC]$-module $M$, the
$\IQ$-vector space $H_p(C_\ast(EW_GC) \otimes_{\IZ[W_GC]} M)$ is isomorphic to
$M \otimes_{\IQ[W_GC]} \IQ$ for $p = 0$ and trivial for $p \geq 1$.
\end{proof}

%This completes the proof of Theorem~\ref{the:positive for IQ and IC}.
%\qed

\subsection{}
{\bf Proof of Theorem~\ref{the:detection-result-for-comm.-rings-of-char.-zero}.}
The proof is analogous to that of Theorem~\ref{the:main-detection-result}, with the exception that
we do not use Proposition~\ref{prop:compute-theta-null} but the following consequences of the hypothesis
on $R$ made in the statement\,: the diagram
\[
\xymatrix{
K_0(\IZ) \ar[r]^-{\dtr}_-{\cong} \ar[d] & \HH_0^{\otimes_{\IZ}} ( \IZ) = \IZ \ar[d]
\\
K_0(R) \ar[r]^-{\dtr}  & \HH_0^{\otimes_{\IZ}}(R)
}
\]
commutes, the upper horizontal map is an isomorphism and both vertical arrows are injective.
The map \eqref{injective-for-z-coefficients} is now defined as the restriction of the upper horizontal arrow
of the diagram appearing in \ref{subsec:general-strategy}, in degree $n$, to the summand for $q = 0$ and
$C = \{e\}$ and then further to the $\IQ$-submodule
$$H_p(BG;\IQ) \cong H_p(BG;\IQ) \otimes_{\IQ} \big( K_0(\IZ) \otimes_{\IZ} \IQ \big) \; \subseteq \;
H_p(BG;\IQ) \otimes_{\IQ} \big( K_0(R) \otimes_{\IZ} \IQ \big)$$
(here, $p=n$). Injectivity of \eqref{injective-for-z-coefficients} is now clear from the general
strategy~\ref{subsec:general-strategy}.
\qed

\subsection{}
{\bf Proof of Theorem~\ref{the:optimal-Z-coeff}.}
For the given $n\geq 0$, the diagram
\begin{eqnarray*}
\qquad\qquad\xymatrix{
H_n^G \big( \EGF{G}{\Fin};\JK^{\minusinfinity} S^{-1} \calo_F \big) \ar[r] \ar[d]&
H_n^G \big( \EGF{G}{\Fin};\JK^{\minusinfinity} F \big) \ar[d]^-{\cong}
\\
H_n^G \big( \EGF{G}{\VCyc};\JK^{\minusinfinity} S^{-1} \calo_F \big) \ar[r] \ar[d]_-{\assembly}^-{\cong_{\IQ}}  &
H_n^G \big( \EGF{G}{\VCyc};\JK^{\minusinfinity} F \big) \ar[d]^-{\assembly}
\\
\hspace*{-7.5em}K_n \big( S^{-1} \calo_F[G] \big) \cong H_n^G(\pt;\JK^{\minusinfinity} S^{-1} \calo_F) \ar[r]
\ar[d]_{\dtr} & H_n^G(\pt;\JK^{\minusinfinity} F) \cong K_n(FG)\hspace*{-4.5em} \ar[d]^{\dtr}
\\
\HH^{\otimes_{\IZ}}_n \big( S^{-1} \calo_F[G] \big) \ar[r]^-{\cong_{\IQ}} &
\HH^{\otimes_{\IZ}}_n(FG)
}
\end{eqnarray*}
commutes, where the vertical maps are induced by the, up to $G$-homotopy, unique $G$-maps or are
given by the trace maps; the horizontal arrows are induced by the inclusion of rings $S^{-1} \calo_F \subseteq F$.
Some explanations are in order for the indicated integral respectively rational isomorphisms.

\smallskip

For every ring $R$, there are isomorphisms
$\HH^{\otimes_{\IZ}}_{\ast} ( R ) \otimes_{\IZ} \IQ \cong
\HH^{\otimes_{\IZ}}_{\ast} ( R \otimes_{\IZ} \IQ )$ and
$\HC^{\otimes_{\IZ}}_{\ast} ( R ) \otimes_{\IZ} \IQ \cong
\HC^{\otimes_{\IZ}}_{\ast} ( R \otimes_{\IZ} \IQ )$, because
$\CN_{\bullet}^{\otimes_{\IZ}} ( R \otimes_{\IZ} \IQ ) \cong
\CN_{\bullet}^{\otimes_{\IZ}} (R) \otimes_{\IZ} \IQ$ and because
the functor $(-)\otimes_{\IZ} \IQ$ commutes with homology and with
$\Tot^{\oplus}$. (For the notation, see
Subsections \ref{subsec:nerves} and \ref{subsec:simpl-dold-kan} below.)
Here, we use that for the total complex occurring in the
definition of cyclic homology it does not matter whether one takes
$\Tot^{\oplus}$ or $\Tot^{\prod}$. Note that a corresponding
statement is false for $\HP$ and $\HN$. Hence the bottom horizontal arrow in the diagram above
is rationally bijective since $S^{-1} \calo_F \otimes_{\IZ} \IQ \cong F$.

\smallskip

The  middle left vertical arrow is rationally  bijective,
since we assume that the Farrell-Jones Conjecture holds rationally for $K_n(S^{-1} \calo_F[G])$.

\smallskip

Since $F$ is a regular ring and contains $\IQ$, the top right vertical arrow is an isomorphism
by \cite[Proposition~2.14]{Lueck-Reich(2003b)}, see also Subsection~\ref{subsec:general-strategy}.

\smallskip

Bartels \cite{Bartels(2003d)} has constructed, for every ring $R$ and every $m\in\IZ$, a retraction
$$r(R)_m \colon H_m^G \big( \EGF{G}{\VCyc};\JK^{\minusinfinity} R \big) \; \to \; H_m^G \big( \EGF{G}{\Fin};
\JK^{\minusinfinity} R \big)$$
of the canonical map $H_m^G(\EGF{G}{\Fin};\JK^{\minusinfinity} R) \to H_m^G(\EGF{G}{\VCyc};\JK^{\minusinfinity} R)$,
which is natural in~$R$. We obtain a decomposition, natural in $R$,
$$H_m^G \big( \EGF{G}{\VCyc};\JK^{\minusinfinity}R \big) \; \cong \;
H_m^G \big( \EGF{G}{\Fin};\JK^{\minusinfinity}R \big) \oplus \ker \big( r(R)_m \big)\,.$$
Therefore, we conclude from the commutative diagram above that the composition
\[ %\label{composition-vcyc-to-HC}
H_n^G \big( \EGF{G}{\VCyc};\JK^{\minusinfinity} S^{-1} \calo_F \big) \xrightarrow{\cong_{\IQ}}
H_n^G ( \pt;\JK^{\minusinfinity} S^{-1} \calo_F ) \xrightarrow{\dtr} \HH^{\otimes_{\IZ}}_n \big( S^{-1} \calo_F[G] \big)\,,
\]
after tensoring with $\IQ$, contains $\ker(r(S^{-1} \calo_F)_n) \otimes_\IZ \IQ$ in its kernel,
because $\ker(r(F)_n)=0$. So, to study injectivity properties of the Dennis trace map
we can focus attention on the composition
\begin{multline*}
H_n^G \big( \EGF{G}{\Fin};\JK^{\minusinfinity} S^{-1} \calo_F \big) \otimes_\IZ \IQ \into
H_n^G \big( \EGF{G}{\VCyc};\JK^{\minusinfinity} S^{-1} \calo_F \big) \otimes_\IZ \IQ
\\
\xrightarrow{\cong} H_n^G ( \pt;\JK^{\minusinfinity} S^{-1} \calo_F ) \otimes_\IZ \IQ
\xrightarrow{\dtr} \HH^{\otimes_{\IZ}}_n \big( S^{-1} \calo_F[G] \big) \otimes_\IZ \IQ\,.
\end{multline*}

\smallskip

By naturality of the bottom isomorphism in Remark~\ref{remark-nonconn-conn}, there is
a commutative diagram
\begin{eqnarray*}
\xymatrix@C=-30mm{
\smash[b]{{\displaystyle {\bigoplus_{\substack{p,q \in \IZ \\p+q = n}}}\; \displaystyle
\bigoplus_{(C) \in (\FCyc)}}}
H_p(BZ_GC;\IQ) \otimes_{\IQ[W_GC]} \theta_C \big(  K_q (S^{-1} \calo_F [C]) \otimes_{\IZ} \IQ \big)
\ar[dr]^-{\ch^G}_-{\cong}
\ar[dd]
& \\
& H_n^G \big( \EGF{G}{\Fin};\JK^{\minusinfinity} S^{-1} \calo_F \big) \otimes_{\IZ} \IQ
 \ar[dd]
\\
\smash[b]{{\displaystyle
{\bigoplus_{\substack{p,q \in \IZ \\p+q = n}}}\; \displaystyle \bigoplus_{(C) \in (\FCyc )}}}
H_p(BZ_GC;\IQ) \otimes_{\IQ[W_GC]} \theta_C \big(  K_q (FC)  \otimes_{\IZ} \IQ \big)
\ar[dr]^-{\ch^G}_-{\cong}
& \\
& H_n^G \big( \EGF{G}{\Fin};\JK^{\minusinfinity} F \big) \otimes_{\IZ} \IQ.
}
\end{eqnarray*}
Now, consider the composition
\begin{multline} \label{composition-oplus-to-HC}
\bigoplus_{\substack{p,q \in \IZ \\p+q = n}}\;
\bigoplus_{(C) \in (\FCyc)}
H_p(BZ_GC;\IQ) \otimes_{\IQ[W_GC]} \theta_C \big(  K_q (S^{-1} \calo_F [C]) \otimes_{\IZ} \IQ \big)
\\[.6em]
\qquad \stackover{\cong\,}{\xrightarrow{\ch^G}}
H_n^G \big( \EGF{G}{\Fin};\JK^{\minusinfinity} S^{-1} \calo_F \big) \otimes_\IZ \IQ
\into H_n^G \big( \EGF{G}{\VCyc};\JK^{\minusinfinity} S^{-1} \calo_F \big) \otimes_\IZ \IQ
\\[1.4em]
\xrightarrow{\cong}
H_n^G ( \pt;\JK^{\minusinfinity} S^{-1} \calo_F ) \otimes_\IZ \IQ
\xrightarrow{\dtr} \HH^{\otimes_{\IZ}}_n \big( S^{-1} \calo_F[G] \big) \otimes_\IZ \IQ\,.
\end{multline}
By the previous two diagrams, the composition \eqref{composition-oplus-to-HC} takes each of the direct summands
for $q \leq -1$ to zero, since $K_q(FC) = 0$ for $q \leq -1$ (the ring $FC$ being regular).

Combining the commutativity of the diagrams occurring in Theorems \ref{the:Chern} and Remark~\ref{remark-nonconn-conn}
(for $R=S^{-1} \calo_F$), we deduce that the composition \eqref{composition-oplus-to-HC} restricted to a
direct summand with $p,q\geq 0$ and with $C$ arbitrary factorizes through the $\IQ$-vector space
\[
H_p(BZ_GC;\IQ) \otimes_{\IQ[W_GC]} \theta_C \big(  \HH^{\otimes_\IZ}_q (S^{-1} \calo_F [C]) \otimes_{\IZ} \IQ \big)\,.
\]
Using the isomorphism $\HH^{\otimes_{\IZ}}_\ast ( S^{-1} \calo_F[G] ) \otimes_\IZ \IQ \cong \HH^{\otimes_{\IZ}}_\ast (FG)
\otimes_\IZ \IQ$, from the vanishing result stated as Proposition~\ref{prop:compute-higher-HH-HC} below, we conclude that
the composition \eqref{composition-oplus-to-HC} vanishes on all summands with $q \geq 1$.

\smallskip

Finally, Proposition~\ref{prop:comp-Z} below implies
that the composition \eqref{composition-oplus-to-HC} vanishes on all summands with $q = 0$ and $C \neq \{e\}$, and is
injective on the summand for $q = 0$ and $C = \{e\}$. But the restriction of the composition
\eqref{composition-oplus-to-HC} to the summand with $q = 0$ and $C = \{e\}$ is precisely the
composition of the injective map \eqref{injective-for-z-coefficients} with the Dennis trace, simply
by Remark~\ref{remark-nonconn-conn} and by construction of the map \eqref{injective-for-z-coefficients}
(see the proof of Theorem~\ref{the:detection-result-for-comm.-rings-of-char.-zero} above). This finishes
the proof of Theorem~\ref{the:optimal-Z-coeff}.
\qed

%%%%%%%%%%%%%%%%%%%%%%%%%%%%%%%%%%%%%%%%%%%%%%%%%%%%%%%%%%%%%%%%%%%%%%%%%%%%%%%%%%%%%%%%%%%%%%%%%%%%%%
%%%%%%%%%%%%%%%%%%%%%%%%%%%%%%%%%%%%%%%%%%%%%%%%%%%%%%%%%%%%%%%%%%%%%%
\section{The trace maps for finite cyclic groups} \label{sec:trace-for-finite-cyclic}

In this section, for a finite cyclic group $C$, a coefficient $k$-algebra $R$, and $q\geq 0$,
we investigate the trace map
\begin{eqnarray} \label{eq:trace-map-with-thetaC}
\theta_C \big( K_q( RC ) \otimes_{\IZ} \IQ \big)  \to \theta_C \big( \HH_q^{\otimes_k} ( RC ) \otimes_{\IZ} \IQ \big)
\end{eqnarray}
and its variants using cyclic, periodic cyclic and negative cyclic homology.
All results concerning the map \eqref{eq:trace-map-with-thetaC} with $q > 0$ will in fact be vanishing results
stating that the map is the zero map.

\begin{remark}
Note that for a commutative ring $k$ and every $k$-algebra $R$,
the canonical maps
\[
\xymatrix{
\HH_0^{\otimes_{\IZ}}(R) \ar[d]_{\cong} \ar[r]^-{\cong} &  \HH_0^{\otimes_k}(R) \ar[d]^{\cong} \\
\HC_0^{\otimes_{\IZ}} (R) \ar[r]^-{\cong} & \HC_0^{\otimes_k} (R)
         }
\]
are all isomorphisms, because all four groups can be identified with $R / [R , R]$.
The following results about $HH_0^{\otimes_{\IZ}}$ hence also apply to other ground rings and to cyclic homology.
\end{remark}

\begin{proposition} \label{prop:compute-theta-null}
Let $G$ be a finite group. Suppose that the ring $R$ is from the following list\,:
\begin{enumerate}
\item \label{injective-i}
a finite dimensional semi\-simple algebra $R$ over a field $F$ of characteristic zero;
\item \label{injective-ii}
a commutative complete local domain $R$ of characteristic zero;
\item \label{injective-iii}
a commutative Dedekind domain $R$ whose quotient field $F$ is an algebraic number field and for
which $|G| \in R$ is invertible.
\end{enumerate}
Then the trace map
$K_0 (R G) \to \HH_0^{\otimes_{\IZ}} ( RG )$ is injective in cases \ref{injective-i} and \ref{injective-ii}  and
is rationally injective in case \ref{injective-iii}. This implies in all cases that for a finite cyclic group $C$,
the induced  map,
\[
\theta_C \big( K_0( RC ) \otimes_{\IZ} \IQ \big)  \to
\theta_C \big( \HH_0^{\otimes_{\IZ}}( RC ) \otimes_{\IZ} \IQ \big)
\]
is injective. Moreover, in all cases, except possibly in case \ref{injective-ii}, the $\IQ$-vector space
$\theta_C ( K_0 (RC) \otimes_{\IZ} \IQ )$ is non-trivial.
\end{proposition}

\begin{proof}
\ref{injective-i}
We first prove injectivity of the trace $K_0( RG ) \to \HH_0^{\otimes_{\IZ}}( RG )$.
Since $R$ is semi\-simple and the order of $G$ is invertible in $R$, the ring $RG$ is
semi\-simple as well, see for example Theorem~6.1 in \cite{Lam(1991)}.
Using the Wedderburn-Artin Theorem \cite[Theorem~3.5]{Lam(1991)} and the fact that the trace map is compatible with
finite products of rings and with  Morita isomorphisms
\cite[Theorem~1.2.4 on page~17 and  Theorem~1.2.15 on page~21]{Loday(1992)},
it suffices to show that the trace map
\[
\dtr\colon K_0 ( D ) \to \HH_0^{\otimes_{\IZ}} ( D )
\]
is injective in the case where $D$ is a skew-field which is a finite dimensional algebra over a field $F$ of
characteristic zero. The following diagram commutes, where the vertical maps are given by restriction to $F$\,:
\[
\xymatrix{
K_0 ( D ) \ar[d]_{\res} \ar[r]^-{\dtr} & \HH_0^{\otimes_{\IZ}} ( D ) \ar[d]^{\res} \\
\; K_0 ( F ) \; \ar[r]^-{\dtr} & \HH_0^{\otimes_{\IZ}} ( F )
         }
\]
The left vertical map can be identified with the map $\dim_F(D) \cdot \id \colon \IZ \to \IZ$
and is hence injective. The trace map $K_0(F) \to \HH_0^{\otimes_{\IZ}}(F)$ can be identified with the inclusion
$\IZ \to F$. This proves injectivity of the Dennis trace $K_0( RG )\to \HH_0^{\otimes_{\IZ}}( RG )$.

\smallskip

Let $R$ be a finite dimensional $F$-algebra. Then induction and restriction with respect to the inclusion
$FG \to RG$ induces maps $\ind \colon K_0(FG) \to K_0(RG)$ and
$\res \colon K_0(RG) \to K_0(FG)$ such that $\res \circ \ind = \dim_F(R) \cdot \id$.
Hence the map $\ind \colon K_0(FG)\otimes_{\IZ} \IQ \to K_0(RG)\otimes_{\IZ} \IQ$ is injective. For $G=C$
a finite cyclic group, this restricts to an injective map
$$ \theta_C \big(K_0(FC)\otimes_{\IZ} \IQ\big) \to \theta_C \big(K_0(RC)\otimes_{\IZ} \IQ\big)\,.$$
Since $F$ is a field of characteristic zero there exists a commutative diagram of ring homomorphisms
\[
\xymatrix{
K_0 ( \IQ C )  \otimes_{\IZ} \IQ \ar[d] \ar[r]^-{\cong} & \map( \sub C , \IQ ) \ar[d] \\
K_0 ( F C ) \otimes_{\IZ} F \ar[r]^-{\cong} & \map (  \Gamma_{F,C}\backslash \con C , F ).
         }
\]
Here, the set $\con C$ of conjugacy classes of elements of $C$ identifies with $C$. Set $m=|C|$ and let
$\mu_m\cong\IZ/m\IZ$ be the group of $m$th roots of $1$ in an algebraic closure of $F$. The action of
the Galois group $G( F(\mu_m) | F )$ on $\mu_m$ determines a subgroup $\Gamma_{F,C}$ of $(\IZ/m \IZ)^{\times}
\cong\Aut (\mu_m)$. An element $t \in \Gamma_{F,C}$ operates on $\con C$ by sending
(the conjugacy class of) the element $c$ to $c^t$.
The set of orbits under this action is $ \Gamma_{F,C}\backslash \con C$.
Note that for $F=\IQ$, the group $\Gamma_{\IQ , C}$ is the whole group
$(\IZ / m \IZ)^{\times}$ and $\Gamma_{\IQ , C}\backslash \con C$ can be
identified with $\sub C$, the set of subgroups of $C$.
So, the first line in the diagram is a special case of the second.
The right-hand vertical map is contravariantly induced from the quotient map
$\Gamma_{F,C}\backslash \con C \to \sub C$ and is in particular injective.
The horizontal maps are given by sending a representation to its character. They are isomorphisms by
\cite[II.~\S~12]{Serre(1977)}.
Hence $\theta_C ( K_0 ( \IQ C ) \otimes_{\IZ} \IQ )$ injects in
$\theta_C ( K_0 ( F C ) \otimes_{\IZ} F )$. We have shown in
Lemma~\ref{lem:compute-theta} that
$\theta_C ( K_0 ( \IQ C ) \otimes_{\IZ} \IQ )$ is non-trivial. Hence
$\theta_C (K_0(RC)\otimes_{\IZ} \IQ)$ is non-trivial as well.

\smallskip

\ref{injective-ii}
 According to Theorem~6.1 in \cite{Swan(1960a)}, the left-hand vertical map in the commutative diagram
\[
\xymatrix{
K_0 ( R G ) \ar[r]^-{\dtr} \ar[d] &  \HH_0^{\otimes_{\IZ}} (RG ) \ar[d] \\
\; K_0 ( FG ) \;\ar[r]^-{\dtr} & \HH_0^{\otimes_{\IZ}} ( FG )
         }
\]
is injective. Here $F$ is the quotient field of $R$. The bottom map is injective by \ref{injective-i}.

\smallskip

\ref{injective-iii} Since any Dedekind ring is regular, the ring
$R$  is a regular domain in which the order of $G$ is invertible.
Hence $RG$ and $FG$ are regular, compare \cite[Proof of
Proposition~2.14]{Lueck-Reich(2003b)}. For any regular ring $S$, the
obvious map $K_0(S) \to G_0(S)$, with $G_0(S)$ the Grothendieck group of
finitely generated $S$-modules, is bijective \cite[Corollary~38.51 on page~29]{Curtis-Reiner(1987)}.
Therefore, the map $K_0 (RG) \to K_0 (FG)$ can be identified with the map
\[
G_0 (R G) \to G_0 ( FG )\,.
\]
This map has a finite kernel and is surjective under our assumptions on
$R$ and $F$ \cite[Theorem~38.42 on page~22 and Theorem~39.14 on page~51]{Curtis-Reiner(1987)}.
We infer that $K_0 ( RG ) \to K_0 ( FG )$ is rationally bijective.
Using the corresponding commutative square involving the trace maps we have reduced our claim to the case \ref{injective-i}.
\end{proof}

\begin{proposition} \label{prop:comp-Z}
Let $S^{-1}\calo_F$ be a localization of the ring of integers $\calo_F$ in an algebraic number field $F$.
Then the canonical map
\[
K_0 ( \IZ ) \otimes_{\IZ} {\IQ} \xrightarrow{\cong} K_0( S^{-1} \calo_F)
\otimes_{\IZ} \IQ
\]
is an isomorphism and the trace map
$$\dtr\colon K_0 ( S^{-1} \calo_F )\otimes_{\IZ} \IQ \to \HH_0^{\otimes_{\IZ}} ( S^{-1} \calo_F )\otimes_{\IZ} \IQ$$
is injective. If $C$ is a non-trivial finite cyclic group and no prime divisor of its order $|C|$
is invertible in $S^{-1}\calo_F$, then
\[
\theta_C ( K_0 ( S^{-1} \calo_F C ) \otimes_{\IZ} \IQ )=0.
\]
\end{proposition}

\begin{proof}
According to a result of Swan \cite[Proposition~9.1]{Swan(1960a)}, the canonical map $K_0 ( \IZ ) \otimes_{\IZ} {\IQ}
\to K_0( S^{-1} \calo_F [G]) \otimes_{\IZ} \IQ$ is an isomorphism for a finite group $G$ if no prime divisor of $|G|
\in\calo_F$ occurs in $S$. As a consequence, the Artin defect \eqref{eq:SC-ist-Artin-defekt} of $K_0 ( S^{-1} \calo_F C )
\otimes_{\IZ} \IQ$ (i.e.\ in degree $0$) vanishes. The result now follows from the identification, which will be proved
in Lemma~\ref{lem:module-over-burnside} below, of $\theta_C ( K_0 ( S^{-1} \calo_F C ) \otimes_{\IZ} \IQ )$
with the Artin defect.
\end{proof}

We next collect the results which state that the trace map is the zero map in higher degrees.
Note that all linear trace maps factorize through $\HN^{\otimes_{\IZ}}_\ast$. The following
result implies that they all vanish in positive degrees for suitable rings $R$.

\begin{proposition} \label{prop:compute-higher-HH-HC}
Let $F$ be an algebraic number field and $R$ a finite dimensional semi\-simple $F$-algebra.
Then, for every finite cyclic group $C$ and for every $n\geq 1$, we have
\[
\HH_n^{\otimes_{\IZ}} ( RC) \otimes_{\IZ} \IQ= 0 \quad \mbox{and} \quad
\HN^{\otimes_{\IZ}}_n ( RC ) \otimes_{\IZ} \IQ  = 0\,.
\]
\end{proposition}

\begin{proof}
Analogously to the proof of
Proposition~\ref{prop:compute-theta-null} \ref{injective-i}, one reduces the claim to the case
where the ring $RC$ is replaced by a skew-field $D$ which is a finite dimensional algebra over
an algebraic number field $F$. Let $\overline{F}$ be a \emph{splitting
field} for $D$, i.e.\ a finite field extension $\overline{F}$ of $F$ such
that $\overline{F} \otimes_F D \cong M_n(\overline{F})$, for some $n\geq 1$, see
\cite[Corollary~7.22 on page~155]{Curtis-Reiner(1981)}. Induction
and restriction for $D \subseteq \overline{F} \otimes_{F} D$ yield maps $\ind
\colon K_\ast(D) \to K_\ast(\overline{F} \otimes_{F} D)$ and $\res \colon
K_\ast(\overline{F} \otimes_F D ) \to K_\ast(D)$ such that $\res \circ
\ind = \dim_F(\overline{F}) \cdot \id$. Hence $\ind \colon K_\ast(D)
\to K_\ast(\overline{F} \otimes_F D)$ is rationally injective. The same
procedure applies to Hochschild homology, cyclic, periodic cyclic and negative
cyclic homology, and all these induction and restriction maps are compatible
with the various trace maps. Applying Morita invariance, it thus suffices to prove that
\[
\HH_n^{\otimes_{\IZ}} (\overline{F}) \otimes_{\IZ} \IQ= 0 \quad \mbox{and} \quad
\HN^{\otimes_{\IZ}}_n (\overline{F}) \otimes_{\IZ} \IQ  = 0\,,
\]
for every $n\geq 1$.
For every $\IQ$-algebra $A$, there is obviously an isomorphism $\CN_{\bullet}^{\otimes_{\IZ}} ( A )
\cong \CN_{\bullet}^{\otimes_{\IQ}}( A )$ of cyclic nerves, see Subsection~\ref{subsec:nerves} for
the notation; hence an isomorphism $\HX_{\ast}^{\otimes_{\IZ}} ( A ) \cong \HX_{\ast}^{\otimes_{\IQ}} ( A )$, where
$\HX$ stands for $\HH$, $\HC$, $\HP$ or $\HN$.
So, we may consider $\HX_\ast^{\otimes_{\IQ}}$ in place of $\HX_\ast^{\otimes_{\IZ}}$ in the sequel.

By the Hochschild-Kostant-Rosenberg Theorem, one has
$\HH_{\ast}^{\otimes_{\IQ}} ( \overline{F} ) \cong \Lambda_{ \overline{F}}^{\ast} \Omega^1_{ \overline{F}| \IQ}$, compare
\cite[Theorem~3.4.4 on page~103]{Loday(1992)}. But $\Omega^1_{ \overline{F} | \IQ} = 0$
because $ \overline{F}$ is a finite separable extension of $\IQ$ (\cite[Corollary~16.16]{Eisenbud(1995)});
therefore $\HH_{\ast}^{\otimes_{\IQ}} ( \overline{F} ) \cong  \overline{F}$ and is concentrated in degree $0$.
From the long exact sequence
\[
\ldots \to \HH_n^{\otimes_{\IQ}}( \overline{F} )\to  \HC_n^{\otimes_{\IQ}}( \overline{F} )\xrightarrow{S}  \HC_{n-2}^{\otimes_{\IQ}}
( \overline{F} )
\to \HH_{n-1}^{\otimes_{\IQ}}( \overline{F} )\to \ldots
\]
it follows that $\HC_{\ast}^{\otimes_{\IQ}} (  \overline{F} )$ is isomorphic to $ \overline{F}$ in each even non-negative degree,
and is zero otherwise. Since the periodicity map $S$ is an isomorphism as soon as its target is non-trivial, the periodic cyclic
homology is the inverse limit $\HP^{\otimes_{\IQ}}_n ( \overline{F} )= \lim_k \HC^{\otimes_{\IQ}}_{n + 2k} (  \overline{F} )$
and hence is concentrated in (all) even degrees, with a copy of $\overline{F}$ in each such degree, compare
\cite[5.1.10 on page~163]{Loday(1992)} and also Remark~\ref{rem:HP-and-Z} below. In the long exact sequence
\[
\ldots \to \HN_{n}^{\otimes_{\IQ}} ( \overline{F} )\to \HP_n^{\otimes_{\IQ}} ( \overline{F} ) \xrightarrow{\overline{S}}
\HC_{n-2}^{\otimes_{\IQ}} ( \overline{F} )\to
\HN_{n-1}^{\otimes_{\IQ}} (  \overline{F} ) \to \ldots
\]
compare~\cite[Proposition~5.1.5 on page~160]{Loday(1992)}, the map $\overline{S}$ is
then an isomorphism whenever its target is non-trivial. It follows
that $\HN^{\otimes_{\IQ}}_\ast ( \overline{F} )$ is concentrated in
non-positive even degrees (with a copy of $\overline{F}$ in each such degree).
\end{proof}

\begin{remark} \label{rem:HP-and-Z}
We could not decide the question whether for an odd $n \geq 1$ and a finite cyclic group $C$, the map
$\theta_C ( K_n ( \IZ C ) \otimes_{\IZ} \IQ ) \to \theta_C ( \HN_n^{\otimes_{\IZ}} ( \IZ C ) \otimes_{\IZ} \IQ)$,
or the corresponding map to periodic cyclic homology, is non-trivial. The calculations in \cite{Kassel(1989)}
and \cite{Cortinas-Guccione-Villamayor(1989)} show that a finer analysis of the trace map is needed in order
to settle the problem. The difficulty is that the $\lim^1$-terms in the computation of $\HP$ out of $\HC$
might contribute to non-torsion elements in odd positive degrees.
\end{remark}

%%%%%%%%%%%%%%%%%%%%%%%%%%%%%%%%%%%%%%%%%%%%%%%%%%%%%%%%%%%%%%%%%%%%%%%%%%%%%%%%%%%%%
%%%%%%%%%%%%%%%%%%%%%%%%%%%%%%%%%%%%%%%%%%%%%%%%%%%%%%%%%%%%%%%%%%%%%%%%%%%%%%%%%%%%%
\section{Notation and generalities}\label{sec:Notation}

%%%%%%%%%%%%%%%%%%%%%%%%%%%%%%%%%%%%%%%%%%%%%%%%%%%%%%%%%%%%%%%%%%%%%%%%%%%%%%%%%%%%%%%%%%%%%%%%%%%%%%
\subsection{Categories and $k$-linear categories} \label{subsec:cat-and-linear-cat}

Let $k$ be a commutative ring. A \emph{$k$-linear category} is a small category which is enriched over
$k$-modules, i.e.\ each morphism set $\hom_{\cala}(c,d)$, with $c,d\in \obj \cala$, has the structure
of a $k$-module, composition of morphisms is bilinear and satisfies the usual associativity axiom; moreover,
there are unit maps $k \to \hom_{\cala} (c , c)$, for every object $c$, satisfying a unit axiom.
Compare \cite[I.8 on page~27, VII.7 on page~181]{MacLane(1971)}. Let $R$ be a $k$-algebra.
For any small category~$\calc$, we can form the associated $k$-linear category $R\calc$.
It has the same objects as $\calc$ and the morphism $k$-modules are obtained as
the free $R$-module over the morphism sets of~$\calc$, i.e.\
\[
\hom_{R \calc} (c,d) = R \mor_{\calc} (c,d)\,.
\]
In fact, this yields a functor $R(-)$  from small categories to $k$-linear categories. Given a $k$-linear category $\cala$,
we denote by $\cala_{\oplus}$ the $k$-linear category whose objects are finite sequences of objects of $\cala$, and
whose morphisms are ``matrices'' of morphisms in $\cala$ with the obvious ``matrix product'' as composition.
Concatenation of sequences yields a sum denoted by ``\,$\oplus$\,'' and we hence obtain, functorially, a
\emph{$k$-linear category with finite sums}, compare \cite[VIII.2 Exercise~6 on page~194]{MacLane(1971)}.
If we consider a $k$-algebra $R$ as a $k$-linear category with one object
then $R_{\oplus}$ is a small model
for the category of finitely generated free left $R$-modules.

%%%%%%%%%%%%%%%%%%%%%%%%%%%%%%%%%%%%%%%%%%%%%%%%%%%%%%%%%%%%%%%%%%%%%%%%%%%%%%%%%%%%%
\subsection{Nerves and cyclic nerves} \label{subsec:nerves}

Let $\calc$ be a small category and let $\cala$ be a $k$-linear category. The \emph{cyclic nerve} of $\calc$
and the \emph{$k$-linear cyclic nerve} of $\cala$ are respectively denoted by
\[
\CN_{\bullet} \calc \quad \mbox{and} \quad \CN_{\bullet}^{\otimes_k} \cala\,.
\]
Depending on the context, they are considered as a cyclic set or as a simplicial set, respectively as a cyclic $k$-module
or as simplicial $k$-module. Recall that by definition, we have
\begin{eqnarray*}
\CN_{[q]} \calc & = &
\coprod_{c_0 , c_1 , \ldots ,c_q \in  \obj \calc} \mor_{\calc} (c_1, c_0) \times \dots \times \mor_{\calc} ( c_q , c_{q-1} )
\times \mor_{\calc}( c_0 , c_q ) ,  \\
\CN_{[q]}^{\otimes_k} \cala & = &
\bigoplus_{c_0 , c_1 , \ldots ,c_q \in  \obj \cala} \hom_{\cala} (c_1, c_0) \otimes_k \dots \otimes_k \hom_{\cala} ( c_q , c_{q-1} )
\otimes_k  \hom_{\cala}( c_0 , c_q ) .  \\
\end{eqnarray*}
The simplicial and cyclic structure maps are induced by composition, insertion of identities and cyclic permutations of
morphisms. For more details, see \cite[2.3]{Waldhausen(1979)}, \cite{Goodwillie(1985)} and \cite{Dundas-McCarthy(1994)}.
The (ordinary) nerve of a small category $\calc$ will always  be considered as a simplicial category and denoted by
$\caln_{\bullet}\calc$. We will write $\obj \caln_{\bullet} \calc$ for the underlying simplicial set of objects.

%%%%%%%%%%%%%%%%%%%%%%%%%%%%%%%%%%%%%%%%%%%%%%%%%%%%%%%%%%%%%%%%%%%%%%%%%%%%%%%%%%%%%%%%%%%%%%%%%%%%%%
\subsection{Simplicial abelian groups and chain complexes} \label{subsec:simpl-dold-kan}

If we are given a simplicial abelian group $M_{\bullet}$, we
denote by $\DK_{\ast}(M_{\bullet})$ the associated normalized
chain complex. For a chain complex of abelian groups $C_{\ast}$
which is concentrated in non-negative degrees, we denote by
$\DK_{\bullet}(C_{\ast})$ the simplicial abelian group that is
associated to it under the Dold-Kan correspondence. For details
see \cite[Section~8.4]{Weibel(1994)}. In particular, recall that
there are  natural isomorphism
$\DK_{\bullet}(\DK_{\ast}(M_{\bullet})) \cong M_{\bullet}$ and
$\DK_{\ast}(\DK_{\bullet}(C_{\ast})) \cong C_{\ast}$.

The \emph{good truncation} $\tau_{\scriptscriptstyle\geq 0} C_{\ast}$ of a chain complex ¤$C_{\ast}$
is defined as the non-negative chain complex which coincides with $C_{\ast}$ in
strictly positive degrees, has the $0$-cycles $Z_0(C_{\ast})$ in degree $0$, and only trivial
modules in negative degrees. Given a bicomplex $C_{\ast \ast}$, we denote by
\[
\Tot^{\oplus} C_{\ast \ast} \quad  \mbox{and} \quad
\Tot^{\prod} C_{\ast \ast}
\]
the total complexes formed using, respectively, the direct sum or the direct product, compare
\cite[1.2.6 on page~8]{Weibel(1994)}.

%%%%%%%%%%%%%%%%%%%%%%%%%%%%%%%%%%%%%%%%%%%%%%%%%%%%%%%%%%%%%%%%%%%%%%%%%%%%%%%%%%%%%%%%%%%%%%%%%%%%%%
\subsection{Spectra, $\Gamma$-spaces and Eilenberg-Mc\,Lane spectra} \label{subsec:Gamma-spaces}

For us, a \emph{spectrum} consists of a sequence $\JE$ of pointed spaces $E_n$, with $n\geq 0$, together with pointed
maps $s_n \colon S^1 \sma E_n \to E_{n+1}$. We do \emph{not} require that the adjoints  $\sigma_n \colon E_n \to \Omega E_{n+1}$ of these maps are homotopy equivalences.
A map of spectra $\Jf \colon \JE \to \JE^{\prime}$ consists of a sequence of maps $f_n \colon E_n \to E_n^{\prime}$ such that
$f_{n+1} \circ s_n  = s_n^{\prime} \circ \id_{S^1} \sma f_n$. One defines in the usual way the homotopy groups
as $\pi_{n}(\JE) = \colim_k \pi_{n+k}( E_k)$, with $n\in\IZ$. The spectrum $\JE$ is \emph{connective} if $\pi_{n}(\JE) = 0$
for all $n<0$. A map of spectra is called a \emph{stable weak equivalence},
or to be short, an \emph{equivalence}, if it induces an isomorphism on all homotopy groups. A \emph{spectrum of simplicial sets}
is defined similarly, using pointed simplicial sets in place
of pointed spaces. Such spectra can be realized and then yield spectra in the sense above.
We denote by $\JS$ the sphere spectrum (as a spectrum of simplicial sets).

Let $\Gamma^{\op}$ denote the small model for the category of finite pointed sets whose objects are $k_{+}=
\{ + , 1, \ldots , k \}$, with $k\geq 0$, and whose morphisms are pointed maps.
A \emph{$\Gamma$-space} $\IE$ is a functor from the category $\Gamma^{\op}$ to the category of
pointed simplicial sets which sends $0_+ = \{ + \}$ to the (simplicial) point. Every $\Gamma$-space $\IE$ can be
extended in an essentially unique way to an endofunctor of the category of pointed simplicial sets
which we again denote by $\IE$. By evaluation on the simplicial spheres, a $\Gamma$-space $\IE$
gives rise to a spectrum of simplicial sets denoted by $\IE ( \JS ) $. The realization $| \IE ( \JS) |$ is then a
spectrum in the sense defined above. A $\Gamma$-space $\IE$ is called \emph{special} if the map
$\IE(k_+) \to \IE(1_+ ) \times\dots \times \IE (1_+)$  induced by the projections
$p_i \colon  k_+ \to 1_+$, with $i=1, \ldots , k$, is a weak equivalence for every $k$.
Here, $p_i ( j )$ is $1$ if $j=i$, and is $+$ otherwise. For more information
on spectra and $\Gamma$-spaces, we refer to \cite{Bousfield-Friedlander(1978)} and \cite{Lydakis(1999)}.

\smallskip

An important example of a $\Gamma$-space is the \emph{Eilenberg-Mc\,Lane $\Gamma$-space} $\IH M_{\bullet}$
associated to a simplicial abelian group $M_{\bullet}$\,. Its value on the finite pointed set $k_+$ is given by
the simplicial abelian group $\IH M_{\bullet} (k_+) = \widetilde{\IZ}[k_+] \otimes_{\IZ} M_{\bullet}$. Here
$\IZ [S]$ denotes the free abelian group generated by the set $S$, and, if the set $S$ is pointed with $s_0$ as
base-point, $\widetilde{\IZ}[S]=\IZ [S]\big/\IZ[s_{0}]$ is the corresponding reduced group. The spectrum
$\JH M_{\bullet} = | \IH M_{\bullet} ( \JS) |$ is a model for the \emph{Eilenberg-MacLane spectrum} associated to
$M_{\bullet}$\,.
The $\Gamma$-space $\IH M_{\bullet}$ is
\emph{very special} in the sense of \cite[page~98]{Bousfield-Friedlander(1978)} and by \cite[Theorem~4.2]{Bousfield-Friedlander(1978)},
the homotopy groups of the spectrum $\JH M_{\bullet}$ coincide with the (unstable) homotopy groups of
(the realization of) $\IH M_{\bullet}( S^0 )$ and hence of $M_{\bullet}$, and consequently with the homology groups of the
associated chain complex $\DK_{\ast}(M_{\bullet})$. So we have natural isomorphisms
\begin{eqnarray} \label{homotopy-HM}
\pi_\ast(\JH M_{\bullet}) \cong \pi_\ast\big(\big|\IH M_{\bullet}( S^0 )\big|\big) \cong \pi_\ast\big(|M_{\bullet}|\big)
\cong H_{*}\big(\DK_{\ast}(M_{\bullet})\big)\,.
\end{eqnarray}

%%%%%%%%%%%%%%%%%%%%%%%%%%%%%%%%%%%%%%%%%%%%%%%%%%%%%%%%%%%%%%%%%%%%%%%%%%%%%%%%%%%%%%%%%%%%%%%%%%%%%%
\subsection{Cyclic, periodic cyclic and negative cyclic homology} \label{define-bicomplexes}

If $Z_{\bullet}$ is a cyclic object in the category of abelian groups, then we denote by $B_{\ast\ast}(Z_{\bullet})$,
$B_{\ast\ast}^{\per}( Z_{\bullet} )$ and $B_{\ast\ast}^{-}( Z_{\bullet} )$ the \emph{cyclic}, \emph{periodic cyclic} and
\emph{negative cyclic bicomplexes}, see~\cite[pages 161--162]{Loday(1992)}. For the good truncations
of the associated total complexes, we write
\begin{eqnarray*}
C_{\ast}^{\HC} (Z_{\bullet}) & = & \Tot^{\prod} B_{\ast\ast}(Z_{\bullet}) \\
C_{\ast}^{\HP} (Z_{\bullet}) & = & \tau_{\scriptscriptstyle\geq 0} \Tot^{\prod} B_{\ast\ast}^{\per}(Z_{\bullet}) \\
C_{\ast}^{\HN} (Z_{\bullet}) & = & \tau_{\scriptscriptstyle\geq 0} \Tot^{\prod} B_{\ast\ast}^{-}(Z_{\bullet})\,.
\end{eqnarray*}
In order to have a uniform notation, it is also convenient to write
\begin{eqnarray*}
C_{\ast}^{\HH} (Z_{\bullet}) & = &  \DK_{\ast} ( Z_{\bullet} ).
\end{eqnarray*}
There is a commutative
diagram of chain complexes
\begin{eqnarray} \label{dia-square-nphc}
\vcenter{\xymatrix{
C_{\ast}^{\HN} ( Z_{\bullet} ) \ar[r] \ar[d]_{{\rm h}_\ast} & C_{\ast}^{\HP} ( Z_{\bullet} ) \ar[d] \\
C_{\ast}^{\HH} ( Z_{\bullet} ) \ar[r] & C_{\ast}^{\HC} ( Z_{\bullet} )
}}
\end{eqnarray}
where the horizontal arrows are induced by inclusions of sub-bicomplexes and the vertical arrows by
projections onto quotient bicomplexes. Let $k$ be a commutative ring. If $Z_{\bullet}$ is the $k$-linear
cyclic nerve $\CN_{\bullet}^{\otimes_k} ( \cala )$ of a $k$-linear category $\cala$, we abbreviate
\begin{eqnarray*}
C_{\ast}^{\HX^{\otimes_{\!k}}}( \cala ) & = &  C_{\ast}^{\HX} \big( \CN_{\bullet}^{\otimes_k} ( \cala ) \big)\,.
\end{eqnarray*}
Here $\HX$ stands for $\HH$, $\HC$, $\HP$ or $\HN$.
The corresponding simplicial abelian group and the corresponding Eilenberg-Mc\,Lane spectrum will be denoted
\begin{eqnarray*}
C_{\bullet}^{\HX^{\otimes_{\!k}}} ( \cala ) & = & \DK_{\bullet} \big( C_{\ast}^{\HX^{\otimes_{\!k}}} (\cala) \big) \\
\JH \JX^{\otimes_k} ( \cala ) & = & \JH C_{\bullet}^{\HX^{\otimes_{\!k}}} (\cala)\,.
\end{eqnarray*}
In particular we have the map $\Jh \colon \JH \JN^{\otimes_k} ( \cala ) \to \JH \JH^{\otimes_k} ( \cala ) $
induced from the map $h_{\ast}$ in \eqref{dia-square-nphc}.
If $R$ is a $k$-algebra we can consider it as a $k$-linear category with one object. Then the
homology groups of $C_{\ast}^{\HX^{\otimes_{\!k}}} ( R )$ as defined above coincide in non-negative degrees with the groups
$\HX^{\otimes_k}_\ast(R)$ that appear in the literature, for instance in \cite{Loday(1992)}. Often negative
cyclic homology $\HN_\ast^{\otimes_k}(R)$ is denoted by $\HC_\ast^-(R)$ or $\HC_\ast^-(R|k)$ in the literature.

%%%%%%%%%%%%%%%%%%%%%%%%%%%%%%%%%%%%%%%%%%%%%%%%%%%%%%%%%%%%%%%%%%%%%%%%%%%%%%%%%%%%%
%%%%%%%%%%%%%%%%%%%%%%%%%%%%%%%%%%%%%%%%%%%%%%%%%%%%%%%%%%%%%%%%%%%%%%%%%%%%%%%%%%%%%
\section{The trace maps} \label{sec:trace-maps}

The aim of this section is to produce the diagram \eqref{spec-level-trace-maps}, i.e.\ the trace maps as maps of
$\Or G$-spectra. We will concentrate on the part of the diagram involving $K$-theory, Hochschild homology and negative
cyclic homology. The remaining arrows are obtained by straightforward modifications.

%%%%%%%%%%%%%%%%%%%%%%%%%%%%%%%%%%%%%%%%%%%%%%%%%%%%%%%%%%%%%%%%%%%%%%%%%%%%%%%%%%%%%
\subsection{The trace maps for additive categories} \label{subsec-trace-additive-categories}

We now review the construction of $K$-theory for additive categories, and of the map $\h$
and the trace maps $\ntr$ and $\dtr$ for $k$-linear categories with finite sums, following the
ideas of \cite{McCarthy(1994)}, \cite{Dundas-McCarthy(1996)} and \cite{Dundas(2000)}.

The following commutative diagram is natural in the $k$-linear category $\cala$\,:
{\small\begin{eqnarray} \label{babytrace}
\hspace*{-4.5em}\vcenter{\mbox{\phantom{xxxxxxxx}}\xymatrix@C=1.7em{
& \qquad\quad C_{[0]}^{\HN^{\otimes_{\!k}}} ( \cala ) = \DK_{0}( C_{\ast}^{\HN^{\otimes_{\!k}}} ( \cala ))
\ar@<-2.5em>[d]^{{\rm h}_0} \ar[r] &  \DK_{\bullet} ( C_{\ast}^{\HN^{\otimes_{\!k}}}( \cala )) = C_{\bullet}^{\HN^{\otimes_{\!k}}}
 ( \cala) \ar@<-2.8em>[d]^{{\rm h}_\bullet=\DK_\bullet({\rm h}_\ast)} \\
\obj \cala \ar[r]^-{\dtr_0} \ar@(dl,u)[ur]^(.46){\ntr_0} & C_{[0]}^{\HH^{\otimes_{\!k}}}( \cala ) = \CN_{0}^{\otimes_k}
\cala \ar[r] & \DK_{\bullet} ( C^{\HH^{\otimes_{\!k}}}_{\ast}( \cala)) \cong \CN_{\bullet}^{\otimes_k}( \cala ).
         }}
\end{eqnarray}}
\hspace*{-.4em}Here, the lower horizontal map $\dtr_0$ is given by sending an object to the corresponding identity
morphism. The lift $\ntr_0$ of this map is explicitly described on page~286 in \cite{McCarthy(1994)}. The remaining
horizontal maps are just the inclusions of the zero simplices. The vertical maps are induced by the map $\h_\ast$ in
diagram \eqref{dia-square-nphc}. The isomorphism in the bottom right corner is a special case
of the natural isomorphism $\DK_{\bullet} ( \DK_{\ast} (M_{\bullet})) \cong M_{\bullet}$, compare Subsection~\ref{subsec:simpl-dold-kan}.
It will be considered as an identification in the following.

The model for the trace maps, for a given $k$-linear category with finite sums $\cala$
will be obtained by replacing
$\cala$ in the diagram above by a suitable simplicial $k$-linear $\Gamma$-category.
On the $K$-theory side, we will use the fact that $\cala$ has finite sums; on the Hochschild
side, we will use the $k$-linear structure.

Let $\cala$ be a small category with finite sums.
We can then apply the Segal construction
which yields a $\Gamma$-category $\Seg \cala$, that is, a functor from $\Gamma^{\op}$ to the
category of small categories, compare~\cite[Definition~3.2]{Dundas(2000)}
and \cite[Section~2]{Segal(1974)}.

Recall that we consider the nerve of a category as a simplicial category.
Let
$\caln_{\bullet}^{\iso}\cala$ be the simplicial subcategory of $\caln_{\bullet} \cala$ for which the objects in
$\caln_{[q]}^{\iso} \cala$ are $q$-tuples of composable isomorphisms in $\cala$, whereas there is no restriction on
the morphisms. Observe that $\obj \caln_{\bullet}^{\iso} \calc = \obj \caln_{\bullet} \iso \calc$, where $\iso \calc$
stands for the subcategory of isomorphisms.

The \emph{connective $K$-theory spectrum}
$\JK( \cala )$ can now be defined as the spectrum associated to the $\Gamma$-space $\obj
\caln_{\bullet}^{\iso} \Seg \cala$, that is,
\begin{eqnarray} \label{def-K-theory}
\JK ( \cala ) & = & \big| ( \obj \caln_{\bullet}^{\iso}  \Seg \cala ) ( \JS ) \big|\,.
\end{eqnarray}
For a comparison to other definitions of $K$-theory, see \cite[Section~1.8]{Waldhausen(1985)}.

We proceed to discuss the trace maps. Let $\cala$ be a $k$-linear category with finite sums. Recall that $\Delta$
is the category of finite ordered sets $[n]=\{0 \leq 1 \leq \ldots \leq  n\}$, with $n\geq 0$, and monotone
maps as morphisms. Observe that $\caln_{\bullet}^{\iso}  \Seg \cala$ is a functor from $\Delta^{\op} \times \Gamma^{\op}$
to $k$-linear categories and it hence makes sense to apply the cyclic nerve constructions.
Since the diagram \eqref{babytrace} is natural in $\cala$
we obtain maps of simplicial $\Gamma$-spaces (alias natural transformations of functors from
$\Delta^{\op} \times \Delta^{\op} \times \Gamma^{\op}$ to the category of pointed sets)
\begin{eqnarray} \label{dia-basic-model-trace-map}
\vcenter{\xymatrix{
& \DK_{\bullet} C_{\ast}^{\HN^{\otimes_{\!k}}} \caln_{\bullet}^{\iso} \Seg \cala \ar[d]^{\h_{\bullet\bullet}} \\
\obj \caln_{\bullet}^{\iso} \Seg \cala
\ar[r]^-{\dtr_{\bullet\bullet}} \ar@(ld,u)[ur]^{\ntr_{\bullet\bullet}} &
\CN_{\bullet}^{\otimes_k} \caln_{\bullet}^{\iso} \Seg \cala.
         }}
\end{eqnarray}
Here $\obj \caln_{\bullet}^{\iso}  \Seg \cala$ is constant in one of the simplicial directions. Taking the
diagonal of the two simplicial directions and passing to the
associated spectra yields the model for the trace maps that we will use. It remains to identify the
objects on the right in \eqref{dia-basic-model-trace-map} with our more standard definitions of Hochschild
and negative cyclic homology.

\begin{lemma} \label{zig-zag-lemma}
Let $\cala$ be a $k$-linear category with finite sums.
There is a zigzag of stable weak equivalences, natural in $\cala$,  between
\[
\xymatrix@R=.5em{
\JH \JN^{\otimes_k} (\cala ) = \JH \DK_{\bullet} C_{\ast}^{\HN^{\otimes_{\!k}}} \cala \ar@<2.4em>[dd]^{\Jh} & &
\big| ( \DK_{\bullet} C_{\ast}^{\HN^{\otimes_{\!k}}} \caln^{\iso}_{\bullet} \Seg \cala) ( \JS ) \big|
\ar[dd]^{|\h_{\bullet\bullet}(\JS)|} \\
& \mbox{and} & \\
\JH \JH^{\otimes_k} (\cala ) = \JH \CN_{\bullet}^{\otimes_k} \cala & &
\big| ( \CN_{\bullet}^{\otimes_k} \caln^{\iso}_{\bullet} \Seg \cala) ( \JS ) \big|.
         }
\]
\end{lemma}

\begin{proof}
Consider, for each $q$, the inclusion of the zero simplices
$i\colon \cala = \caln_{[0]}^{\iso} \cala \to \caln_{[q]}^{\iso} \cala$. There is a
left inverse $p$ (forget everything but the zero-th object) and an obvious natural transformation
between $i \circ p$ and the identity which is objectwise an isomorphism. This induces a special
homotopy equivalence (\cite[Definition~2.3.2]{McCarthy(1994)}) and hence in particular a homotopy
equivalence of cyclic nerves $\CN_{\bullet}^{\otimes_k} \cala \xrightarrow{\simeq} \CN_{\bullet}^{\otimes_k}
\caln_{\bullet}^{\iso} \cala$, which passes to a chain homotopy equivalence
on the negative cyclic construction, compare \cite[Proposition~2.4.1]{McCarthy(1994)}. So
we get rid of the $\caln_{\bullet}^{\iso}$'s in the expressions above. The rest now follows by applying
the following lemma to the map
\[
\IH(h_{\bullet})\colon\IH \DK_{\bullet} C_{\ast}^{\HN^{\otimes_{\!k}}} \Seg \cala \to \IH
\CN_{\bullet}^{\otimes_k}  \Seg \cala
\]
of bi-$\Gamma$-spaces, provided we can prove that the source and the target are both special in both variables
(see~\ref{subsec:Gamma-spaces}). Specialness in the Eilenberg-Mc\,Lane-variable is standard and follows immediately
from the definition of the functor $\IH(-)$. Being special in the Segal-variable means in the case of the first
bi-$\Gamma$-space that for every $l_+$ and $k_+$, the following composition is a weak equivalence\,:
\begin{eqnarray*}
\widetilde{\IZ} [l_+] \otimes_{\IZ} \big(
\CN^{\otimes_k}_{\bullet} \Seg \cala (k_+) \big) & \to &
\widetilde{\IZ} [l_+] \otimes_{\IZ} \big(
\CN^{\otimes_k}_{\bullet} ( \cala \times \dots \times \cala ) \big) \\[.2em]
& \to &
\widetilde{\IZ} [l_+] \otimes_{\IZ} \big(
\CN^{\otimes_k}_{\bullet} \cala \times \dots \times \CN^{\otimes_k}_{\bullet} \cala \big) \\[.2em]
& \to &
\widetilde{\IZ} [l_+] \otimes_{\IZ}
\CN^{\otimes_k}_{\bullet} \cala \times \dots \times
\widetilde{\IZ} [l_+] \otimes_{\IZ}
\CN^{\otimes_k}_{\bullet} \cala\,.
\end{eqnarray*}
This is clearly true for the last map.
The Segal construction is designed in such a way that $\Seg \cala ( k_+ ) \to \cala \times \dots \times \cala$
is an equivalence of categories. By \cite[Proposition~2.4.1]{McCarthy(1994)}, this passes to an equivalence
on the cyclic constructions and yields that the first map is an equivalence.
Proposition~2.4.9 in \cite{McCarthy(1994)} deals with the second map. The argument for the second bi-$\Gamma$-space is
analogous.
\end{proof}

\begin{lemma} \label{bi-special-bi-gamma}
Suppose that $(k_+,l_+) \mapsto \IA (k_+ , l_+)$ is a bi-$\Gamma$-space which is special
in both variables, i.e.\ for every fixed $l_+$, the $\Gamma$-space $k_+ \mapsto \IA ( k_+ , l_+ )$
is special, and similarly in the other variable. Then there is a natural zigzag of
stable weak equivalences of spectra of simplicial sets between
\[
\IA ( 1_+ , \JS ) \quad \mbox{and} \quad \IA ( \JS , 1_+ )\,.
\]
\end{lemma}

\begin{proof}
There is a naive definition of a bi-spectrum (of simplicial sets) as a collection $\IE$ of pointed simplicial sets $E_{n,m}$,
with $n \geq 0$ and $m \geq 0$, together with horizontal and vertical pointed structure maps
$\sigma_h \colon E_{n,m} \sma S^1 \to E_{n+1,m}$ and $\sigma_v \colon S^1 \sma   E_{n,m} \to E_{n,m+1}$
satisfying
\[
\sigma_h \circ (\sigma_v \sma S^1) = \sigma_v \circ (S^1 \sma \sigma_h)\,.
\]
After some choice of a poset map $\mu \colon \IN_0 \to\IN_0 \times \IN_0$ satisfying a suitable cofinality condition,
for example $\mu(2n)=(n,n)$ and $\mu(2n+1)=(n+1,n)$, one can form the diagonal spectrum
$\diag_{\mu} \IE$, compare \cite[Section~1.3]{Jardine(1997)}.

A \emph{bi-$\Gamma$-space} $\IA$ is a functor from $\Gamma^{\op}$ to the category of $\Gamma$-spaces, denoted
by $k_+\mapsto\IA(k_+,-)$, and such that $\IA(0_+,l_+)$ is the (simplicial) point for each $l_+$.
Every bi-$\Gamma$-space $\IA$ gives rise to a simplicial bi-spectrum $\IA(\JS , \JS^{\prime})$ in the naive
sense above, with $\IA(\JS , \JS^{\prime})_{n,m}=\IA ( S^n , S^m )$. Here $\JS^{\prime}$ is
just a copy of the simplicial sphere spectrum $\JS$ which we want to distinguish in the notation.

There are maps of bi-$\Gamma$-spaces
\[
\xymatrix{
\IA ( 1_+ , k_+ \sma l_+ ) & k_+ \sma \IA( 1_+ , l_+ ) \ar[l] \ar[r] & \IA ( k_+ , l_+ )
         }
\]
which extend to maps of bi-spectra. We claim that for every pointed simplicial set $Y$, the corresponding maps of simplicial
spectra
\begin{eqnarray} \label{left-and-right}
\xymatrix{
\IA ( 1_+ , Y \sma \JS ) & Y \sma \IA( 1_+ , \JS) \ar[l] \ar[r] & \IA ( Y , \JS )
         }
\end{eqnarray}
are stable weak equivalences. For the first map, this is \cite[Lemma~4.1]{Bousfield-Friedlander(1978)}
and no specialness assumption is needed. For the second map, one argues as follows.
For a pointed simplicial set $X$, the simplicial set $\IA ( 1_+ , X)$ is at least as connected
as $X$ \cite[Proposition~5.20]{Lydakis(1999)}. Now, the composition
\[
k_+ \sma \IA( 1_+ , X ) \to \IA( k_+ , X ) \to \IA ( 1_+ , X ) \times \dots \times \IA
( 1_+ , X )
\]
is the inclusion of a $k$-fold wedge into the corresponding $k$-fold product and hence roughly twice as connected as
$\IA (1_+ , X)$. Since the second map is a weak equivalence (by the assumption that $\IA$ is special in the
second variable), we conclude that the connectivity of the first map grows faster than $n$ for $X=S^n$.
The same statement holds for arbitrary pointed simplicial sets $Y$ in place of $k_+$
by a careful version of the Realization Lemma for bisimplicial sets (realization preserves connectivity --- compare
\cite[Lemma~2.1.1]{Waldhausen(1979)}). So, the second map in \eqref{left-and-right} indeed is a stable weak equivalence,
proving the claim above.

If we now apply the elementary Lemma~1.28 from \cite{Jardine(1997)} (this is a Realization Lemma for bi-spectra), we
obtain a zigzag of weak equivalences of spectra of simplicial sets between
$\IA ( 1_+ , \diag_{\mu} \JS \sma \JS^{\prime} )=\diag_{\mu} \IA ( 1_+ , \JS \sma \JS^{\prime} )$
and $\diag_{\mu}\IA ( \JS , \JS^{\prime} )$. The pointed isomorphism between
$S^0$ and the $0$th simplicial set of the spectrum $\diag_{\mu} \JS \sma \JS^{\prime}$ determines
uniquely a map of spectra $\JS \to \diag_{\mu} \JS \sma \JS^{\prime}$ which clearly is an isomorphism. In
total, we have constructed a zigzag of stable weak equivalences between $\IA( 1_+ , \JS)$
and $\diag_{\mu} \IA(\JS , \JS^{\prime})$. The result now follows by symmetry (using specialness in the first variable).
\end{proof}

Summarizing we have that for a $k$-linear category $\cala$ with finite sums, the model for the
trace maps, at the level of spectra, is given by a  commutative diagram of the following form
{\small\[
\xymatrix@C=1.2em{
& \big| ( \DK_{\bullet} C_{\ast}^{\HN^{\otimes_{\!k}}} \caln^{\iso}_{\bullet} \Seg \cala) ( \JS ) \big|
\ar[d]^{|\h_{\bullet\bullet}(\JS)|}
\ar[r]^-{\simeq} & & \dots \ar[l]_-{\simeq} \ar[r]^-{\simeq} & \JH \JN^{\otimes_k} (\cala ) \ar[d]^{\Jh}  \\
\JK ( \cala ) = \big| ( \obj \caln_{\bullet}^{\iso}  \Seg \cala ) ( \JS ) \big|
\ar[r]^-{\Jdtr} \ar@(ld,u)[ur]^{\Jntr} &
\;\;\; \big| ( \CN_{\bullet}^{\otimes_k} \caln^{\iso}_{\bullet} \Seg \cala) ( \JS ) \big| \;\;\; \ar[r]^-{\simeq} & &
\dots \ar[l]_-{\simeq} \ar[r]^-{\simeq} & \JH \JH^{\otimes_k} (\cala ).
         }
\]}

%%%%%%%%%%%%%%%%%%%%%%%%%%%%%%%%%%%%%%%%%%%%%%%%%%%%%%%%%%%%%%%%%%%%%%%%%%%%%%%%%%%%%
\subsection{The trace maps as maps of spectra over the orbit category} \label{subsec:trace-maps}

We will now define the $\Or G$-spectra representing $K$-theory, Hochschild homology and other
cyclic homology theories, and the trace maps which appear in \eqref{spec-level-trace-maps}.

\smallskip

Given a $G$-set $S$, let $\calg^G(S)$ denote the associated \emph{transport groupoid}, i.e.\ the category
whose objects are the elements of $S$ and where the set of morphisms from $s \in S$ to $t \in S$ is given
by $\mor(s,t)=\{ g \in G \,|\, gs=t \}$. Given  a $k$-algebra $R$ we can compose the functor $\calg^G(?)$
with the functors $R( -)$ and $( - )_{\oplus}$ (compare Subsection~\ref{subsec:cat-and-linear-cat})
to obtain a functor
\begin{eqnarray*}
R \calg^G ( ? )_{\oplus}\colon\Or G \to  k\mbox{-}\Cat_{\oplus}\,,\quad G/H \mapsto R \calg^G ( G/H )_{\oplus}\,,
\end{eqnarray*}
where $k\mbox{-}\Cat_{\oplus}$ denotes the \emph{category of small $k$-linear categories with finite sums},
whose morphisms are $k$-linear functors (and hence respect the sum, compare \cite[VIII.2~Prop.4 on page 193]{MacLane(1971)}). The \emph{idempotent completion} $\Idem \cala$
of a category $\cala$ has as objects the idempotent endomorphisms in $\cala$, i.e.\ morphisms $p \colon  c \to c$
with $p\circ p=p$; a morphism from $p\colon c \to c$ to $q\colon d \to d$ is given
by a morphism $f\colon c \to d$ with $q\circ f= f \circ p$. The idempotent completion of a $k$-linear category
is again $k$-linear. For a small category $\calc$, the idempotent completion of $R \calc_{\oplus}$ is a $k$-linear
category with finite sums. For an arbitrary ring $S$, the category $\Idem S_{\oplus}$
is a small model for the category of finitely generated projective left $S$-modules.

Let $R$ be a $k$-algebra and $H$ a subgroup of $G$. Consider the commutative diagram of $k$-linear categories
\begin{eqnarray} \label{diagram-k-lin-equivalences}
\hspace*{3em}\vcenter{\xymatrix@C=2em{
RH \ar[r] \ar[d] & RH_{\oplus} \ar[r] \ar[d] & \Idem RH_{\oplus} \ar[d] \\
R\calg^G( G/H )  \ar[r] & R\calg^G (G/H)_{\oplus} \ar[r] & \Idem R\calg^G (G/H)_{\oplus} .
         }}
\end{eqnarray}
The vertical functors are all induced from considering $H$ as the full subcategory of $\calg^G( G/H)$
on the object $eH \in G/H= \obj \calg^G(G/H)$. All vertical functors are $k$-linear equivalences
and the two right hand functors are cofinal inclusions into the corresponding  idempotent
completions. Hence it follows from \cite[Proposition~2.4.1 and 2.4.2]{McCarthy(1994)} that all functors
in the diagram above induce equivalences if one applies Hochschild homology or one of the cyclic homology
theories, i.e.\ $\JH \JX^{\otimes_k} ( - )$. Observe that our $K$-theory functor $\JK( - )$ can only be
applied to the four categories on the right (they have finite sums). The two right-hand vertical maps
induce isomorphisms on all higher $K$-groups, however, $K_0$ may differ for a category with finite sums and
its idempotent completion.

\smallskip

Finally, define  $\Or G$-spectra $\JK R(?)$ and $\JH \JX^{\otimes_k}R(?)$ by
\begin{eqnarray}
\JK R (G/H) & = & \JK \Idem R \calg^G( G/H )_{\oplus} \label{def-KR-OrG-spectrum} \\
\JH \JX^{\otimes_k} R (G/H) & = & \JH \JX^{\otimes_k} R \calg^G( G/H )\,. \label{def-HX-OrG-spectrum}
\end{eqnarray}
Here again $\JH \JX$ stands for $\JH \JH$, $\JH \JC$ , $\JH \JP$ or $\JH \JN$.
Compare \eqref{def-K-theory} and the notation introduced in Subsection~\ref{define-bicomplexes}.
The discussion above and the one in Subsection~\ref{define-bicomplexes} verify all the isomorphisms
claimed in \eqref{relate-to-classical}.

\smallskip

Now, apply the construction of \eqref{dia-basic-model-trace-map} in the case where the additive category
$\cala$ is $\Idem R \calg^G( G/H )_{\oplus}$. Using the equivalences (discussed above) induced by the map
$\Idem R \calg^G( G/H )_{\oplus}\leftarrow R\calg^G( G/H )$ and the equivalences appearing in the
diagram at the end of Subsection~\ref{subsec-trace-additive-categories}, we obtain a commutative diagram
of connective $\Or G$-spectra of the shape
\[
\xymatrix@C=1.7em{
& & {}^{\prime}\JH\JN^{\otimes_k}R \ar[d]& \ar[l]_-{\simeq} {}^{\prime\prime}\JH\JN^{\otimes_k}R \ar[d]
\ar[r]^-{\simeq} & & \dots \ar[l]_-{\simeq} \ar[r]^-{\simeq} & \JH \JN^{\otimes_k} R \ar[d]^{\Jh} \\
\JK R \ar[rr]^-{\Jdtr} \ar[urr]^{\Jntr} & & {}^{\prime}\JH\JH^{\otimes_k}R & \ar[l]_-{\simeq} {}^{\prime\prime}\JH
\JH^{\otimes_k}R \ar[r]^-{\simeq} & & \dots \ar[l]_-{\simeq} \ar[r]^-{\simeq} & \JH \JH^{\otimes_k} R,
         }
\]
where all arrows labelled with a ``\,$\simeq$\,'' (in particular all those pointing left) are objectwise
stable weak equivalences.

%%%%%%%%%%%%%%%%%%%%%%%%%%%%%%%%%%%%%%%%%%%%%%%%%%%%%%%%%%%%%%%%%%%%%%%%%%%%%%%%%%%%%%%%%%%%%%%%%%%%%%
%%%%%%%%%%%%%%%%%%%%%%%%%%%%%%%%%%%%%%%%%%%%%%%%%%%%%%%%%%%%%%%%%%%%%%%%%%%%%%%%%%%%%%%%%%%%%%%%%%%%%%
\section{Equivariant homology theories, induction  and Mackey structures} \label{sec:mackey}

A \emph{$G$-homology theory} is a collection of functors $\calh_\ast^{G}(-)=\{\calh_n^{G}(-)\}_{n\in\IZ}$
from the category of (pairs of) $G$-$CW$-complexes to the category of abelian groups, which satisfies the $G$-analogues
of the usual axioms for a generalized homology theory, compare \cite[2.1.4]{Lueck-Reich(2003b)}.

\smallskip

For example, every $\Or G$-spectrum $\JE=\JE(?)$ gives rise to a $G$-homology theory $H_\ast^G ( - ; \JE )$ by setting,
for a $G$-$CW$-complex $X$,
\begin{eqnarray*}
H_\ast^G ( X ; \JE ) & = & \pi_\ast \big( X_+^? \sma_{\Or G} \JE(?) \big)
\end{eqnarray*}
and more generally, for a pair of $G$-$CW$-complexes $(X,A)$,
\begin{eqnarray*}
H_\ast^G ( X , A ; \JE ) & = & \pi_\ast \big( (X_+/A_+)^? \sma_{\Or G} \JE(?) \big)\,.
\end{eqnarray*}
Here, for a $G$-space $Y$, the symbol $Y_+$ denotes the space $Y$ with a disjoint base-point added (viewed as a
$G$-fixpoint), and $Y^?$ stands for the \emph{fixpoint functor} $\map_G( - , Y)$, considered as a contravariant functor
from $\Or G$ to the category of spaces; and $X_+^? \sma_{\Or G} \JE (?)$ is the \emph{balanced smash product} of
a contravariant pointed $\Or G$-space and a covariant $\Or G$-spectrum. It is constructed by applying levelwise the
\emph{balanced smash product}
{\small\begin{multline}  \label{balanced-smash}
Y \sma_{\Or G} Z  =
\xymatrix@C=1.5em{
 \coequ \Big(
\bigvee_{f \in \mor \Or G} Y( t(f) ) \sma Z( s(f) ) \ar@<3pt>[r] \ar@<-3pt>[r] &
\bigvee_{G/H \in \Or G} Y( H ) \sma Z( H ) \Big)
}
\end{multline}}
\hspace*{-.3em}of a contravariant pointed $\Or G$-space $Y(?)$ and a covariant pointed $\Or G$-space $Z(?)$;
here, $s(f)$ stands for the source and $t(f)$ for the target of the morphism $f \in \mor \Or G$,
$\coequ$ is the \emph{coequalizer}, and the two indicated maps are defined by $f^{\ast} \sma \id$
and $\id \sma f_{\ast}$ on the wedge-summand corresponding to $f$. We repeat that $H_{\ast}^G( \pt ; \JE)$ identifies
with $\pi_{\ast}( \JE( G/G ))$. For details, we refer to \cite{Davis-Lueck(1998)} and \cite[Chapter~6]{Lueck-Reich(2003b)}.

\smallskip

For a group homomorphism $\alpha \colon H \to G$ and an $H$-$CW$-complex $X$, let $\ind_{\alpha} X$ be the quotient
of $G \times X$ by the right action of $H$ given by $(g,x)h = (g \alpha (h) , h^{-1} x)$. An \emph{equivariant homology
theory} $\calh^?_\ast=\calh^?_\ast(-)$ consists of a $G$-homology theory for each group $G$ together with natural induction
isomorphisms
\[
\ind_{\alpha} \colon \calh^H_\ast ( X , A ) \xrightarrow{\cong} \calh^G_\ast ( \ind_{\alpha} X , \ind_{\alpha} A )
\]
for each group homomorphism $\alpha \colon H \to G$ and each $H$-$CW$-pair $(X,A)$ such that $\ker \alpha$ acts freely on $X$.
The induction isomorphisms need to verify certain natural axioms, compare \cite[6.1]{Lueck-Reich(2003b)}. We refer to
the collection of induction isomorphisms as an ``\emph{induction structure}''.

\smallskip

Suppose that an $\Or G$-spectrum $\JD(?)$ is a composition of functors $\JD = \JE \circ \calg^G (?)$, where
$\JE \colon \Groupoids \to \Sp$ is a functor from the \emph{category of small groupoids} to the \emph{category of spectra}.
If $\JE$ is a \emph{homotopy functor}, i.e.\ sends equivalences of groupoids to stable weak equivalences of spectra,
then, according to \cite[Proposition~6.10]{Lueck-Reich(2003b)} and \cite{SauerJ(2002)}, there is a `naturally' defined
induction structure for the collection of $G$-homology theories, one for each group $G$, given by $H_\ast^G ( - ; \JE
\circ \calg^G )$. Hence each homotopy functor $\JE \colon \Groupoids \to \Sp$ determines an equivariant homology theory
$H_\ast^? ( - ; \JE \circ \calg^G )$.

\smallskip

Given an equivariant homology theory $\calh_\ast^?(-)$, one can, for each $n\in\IZ$, construct a covariant functor
from $\FGINJ$, i.e.\ \emph{the category of finite groups and injective group homomorphisms}, to $\Ab$, i.e.\ the
\emph{category of abelian groups}, by setting
\begin{eqnarray} \label{mackey-co-a}
M_{\ast} \colon \FGINJ \to \Ab\,, \quad G \mapsto \calh^G_n( \pt )\,;
\end{eqnarray}
for a group monomorphism $\alpha \colon H \into G$, we define $M_{\ast}(\alpha)$ as the composition
\begin{eqnarray} \label{mackey-co-b}
\xymatrix{
M_{\ast} ( H ) = \calh^H_n ( \pt ) \ar[r]^-{\ind_{\alpha}} & \calh_n^G \big( G/ \alpha(H) \big) \ar[r]^-{\calh_n^G (\pr)} &
\calh^G_n ( \pt ) = M_{\ast} (H)\,,
         }
\end{eqnarray}
where $pr$ is the projection onto the point.

\smallskip

A \emph{Mackey functor} $M$ is a pair $(M_{\ast} , M^{\ast})$ consisting of a co- and a contravariant functor $\FGINJ \to \Ab$
which agree on objects, i.e.\ $M_{\ast}( H ) = M^{\ast} ( H )$ (merely denoted by $M(H)$), and satisfy the following axioms.
\begin{enumerate}
\item For an inner automorphism $c_g \colon G \to G$, $h \mapsto g^{-1} h g$ with $g \in G$ one has
$M_{\ast} ( c_g ) = \id \colon M(G) \to M(G)$.
\item
If $f \colon G \xrightarrow{\cong} H$ is an isomorphism, then one has $M_{\ast} (f) \circ M^{\ast} (f) = \id$ and
$M^{\ast} (f) \circ M_{\ast} (f) = \id$.
\item
There is a double coset formula, i.e., for two subgroups $H,K \leq G$, one has
$$
\renewcommand{\arraystretch}{1.5}
\begin{array}{l}
\qquad\quad\;\; M^{\ast} \big( i \colon K \to G \big) \circ M_{\ast} \big( i \colon H \to G \big) \; = \qquad \\[.2em]
\qquad\qquad\quad\;  \; {\displaystyle\sum_{KgH \in K\backslash G / H}} M_{\ast} \big( c_g \colon H \cap g^{-1} K g \to K \big)
\circ M^{\ast} \big( i \colon H \cap g^{-1} K g \to H \big)\,,
\end{array}
$$
where $c_g( h) = g^{-1}hg$ and $i$ in each case denotes the inclusion.
\end{enumerate}

\smallskip

If, for every $n\in\IZ$, the covariant functor $M_{\ast}$ we associated in
\eqref{mackey-co-a} and  \eqref{mackey-co-b} to an equivariant homology theory $\calh^?_\ast(-)$
can be extended to a Mackey functor, then we say that the equivariant homology theory admits a
``\emph{Mackey structure}''.

\smallskip

Let $R$ be a $k$-algebra. We will consider compositions of functors of the form
\[
\xymatrix@C=3em{
\Or G \ar[r]^-{\calg^G( - )} & \Groupoids \ar[r]^-{R ( - )_{\oplus}} &
k\mbox{-}\Cat_{\oplus} \ar[r]^-{\JF} & \Sp\,.
         }
\]
Recall that $k\mbox{-}\Cat_{\oplus}$ denotes the category of small $k$-linear categories with finite sums.

The $\Or G$-spectra we are mainly interested in, namely $\JK R(?)$ and  $\JH \JX^{\otimes_k}R(?)$,
are defined (up to equivalence for the latter) as such a composition with $\JF$ being the composite functor
$\JK\circ\Idem(-)$ for the former, see \eqref{def-K-theory} and \eqref{def-KR-OrG-spectrum}, and being
$\JH \JX^{\otimes_k}(-)$ for the latter, see Subsection~\ref{define-bicomplexes} and  the discussion following
diagram \eqref{diagram-k-lin-equivalences}, and \eqref{def-HX-OrG-spectrum}. It turns out that the non-connective
$K$-theory $\Or G$-spectrum $\JK^{\minusinfinity}R(?)$ of Example~\ref{exa:FJ} is also such a composition.
In that case $\JF$ is the \emph{Pedersen-Weibel functor} (defined on $\Cat_{\oplus}$), compare \cite{Pedersen-Weibel(1985)}.
Up to equivalence a model for the $(-1)$-connective covering map of $\Or G$-spectra $\JK R(?)\to\JK^{\minusinfinity}R(?)$ mentioned in~\ref{exa:FJ} is induced by
a specific natural transformation between the corresponding $\JF$'s.

So, consider a functor $\JF\colon k\mbox{-}\Cat_{\oplus} \to \Sp$. We call $\JF$ a \emph{homotopy functor} if it takes
$k$-linear equivalences to stable weak equivalences of spectra. We call $\JF$ \emph{additive} if for every $k$-linear
functors $f,g \colon \cala \to \calb$ between $k$-linear categories with finite sums,
\begin{eqnarray} \label{the-sum-property}
\pi_{\ast} ( \JF ( f \oplus g ) ) = \pi_{\ast} ( \JF (f) ) + \pi_{\ast} ( \JF (g) )
\end{eqnarray}
holds; here, $f \oplus g \colon \cala \to \calb$ is the composition
\[
\xymatrix{
\cala \ar[r]^-{\diag} & \cala \times \cala \ar[r]^-{f \times g} & \calb \times \calb \ar[r]^-{\oplus} & \calb\,,
         }
\]
where $\diag$ denotes the \emph{diagonal embedding} and $\oplus$ is the sum in $\calb$.

\begin{proposition} \label{prop:homotopy-functor}
Suppose that $\JF\colon k\mbox{-}\Cat_{\oplus} \to \Sp$ is a homotopy functor. Then, the composite functor
$\JF \circ R (- )_{\oplus}$ is a homotopy functor; in particular, it determines an equivariant homology
theory whose underlying $G$-homology theory, for a group~$G$, is given by the $\Or G$-spectrum $\JF R \calg^G
(?)_{\oplus}$\,, that is, by
\begin{eqnarray*}
H_\ast^G \big( X , A ; \JF R \calg^G(?)_{\oplus} \big) & = & \pi_\ast \big( (X_+/A_+)^? \sma_{\Or G} \JF R\calg^G
( ? )_{\oplus} \big)\,.
\end{eqnarray*}
If $\JF$ is additive then this equivariant homology theory admits a Mackey structure.
\end{proposition}

\begin{proof}
The first part is clear. For the second, we need to define the contravariant half of the Mackey functor and verify the axioms.
For a given ring $S$ let  $\calf(S)$ denote the \emph{category of finitely generated free left $S$-modules},
which is of course not a small category. If we consider a group $H$ as a groupoid with one object, then  $RH_{\oplus}$
is a small model for the category of finitely generated free
left $RH$-modules and there is an inclusion functor $i_H \colon RH_{\oplus} \to \calf(RH)$ which is
an equivalence of categories. We choose a functor $p_H \colon \calf(RH) \to RH_{\oplus}$ such that
$p_H \circ i_H \simeq \id$ and $i_H \circ p_H \simeq \id$. Here, $f \simeq g$ indicates that there
exists a natural transformation through isomorphisms.
Given a homomorphism $\alpha \colon H \to G$ between finite groups, there are the usual induction
and restriction functors $\ind_{\alpha} \colon \calf( RH) \to \calf (RG)$ and $\res_{\alpha} \colon
\calf (RG ) \to \calf ( RH )$. For $n\in\IZ$, we define induction and restriction homomorphisms
\[
\ind_{\alpha} \colon \pi_n \big( \JF RH_{\oplus} \big) \to \pi_n \big( \JF RG_{\oplus} \big) \quad \mbox{and} \quad
\res_{\alpha} \colon \pi_n \big( \JF RG_{\oplus} \big) \to \pi_n \big( \JF RH_{\oplus} \big)
\]
as $\ind_{\alpha} = \pi_n ( \JF ( p_G \circ \ind_{\alpha} \circ i_H ))$ and
$\res_{\alpha} = \pi_n ( \JF ( p_H \circ \res_{\alpha} \circ i_G ))$.
Since $f \simeq g$ implies $\pi_n ( \JF ( f) ) = \pi_n ( \JF (g) )$, this does not depend on the
choice of $p_H$ and $p_G$.

\smallskip

Unravelling the definitions, one checks that under the identifications
\[
M(H) = \pi_n \big( \pt_+^? \sma_{\Or G} \JF R \calg^G( ? )_{\oplus} \big) \cong \pi_n \big( \JF R \calg^G ( H/H )_{\oplus} \big)
\cong \pi_n \big( \JF RH_{\oplus} \big)\,,
\]
the induction homomorphism $ M_{\ast} ( \alpha )$ from \eqref{mackey-co-b} coincides with the
induction homomorphism we have just constructed. Using the same identifications, we consider the map
$\res_{\alpha}$ constructed above as a map $M(G) \to M(H)$ and denote it by $M^{\ast}( \alpha )$. The
axioms now follow since each of the remaining equalities corresponds to a well-known natural isomorphism
between functors on categories of finitely generated free left modules; for the third axiom, one uses
\eqref{the-sum-property}, i.e.\ additivity of $\JF$.
\end{proof}

The functors $\JF$ that are responsible for $\JK^{-\infty}(?)$, $\JK(?)$ and $\JH \JX^{\otimes_k} (?)$ are
homotopy invariant and additive.
We hence obtain the corresponding equivariant homology theories with Mackey structures
given, at a group $G$, by $H_\ast^G ( - ; \JK R )$, $H_\ast^G ( - ; \JK^{\minusinfinity} R )$ and by $H_\ast^G ( - ;\JH \JX^{\otimes_k} R )$.
The maps between these theories that are induced from the maps of $\Or G$-spectra that we have discussed above
are  compatible with the
induction and Mackey structures.

%%%%%%%%%%%%%%%%%%%%%%%%%%%%%%%%%%%%%%%%%%%%%%%%%%%%%%%%%%%%%%%%%%%%%%%%%%%%%%%%%%%%%%%%%%%%%%%%%%%%%%
%%%%%%%%%%%%%%%%%%%%%%%%%%%%%%%%%%%%%%%%%%%%%%%%%%%%%%%%%%%%%%%%%%%%%%%%%%%%%%%%%%%%%%%%%%%%%%%%%%%%%%
\section{Evaluating the equivariant Chern character}
\label{sec:eval-chern}
In this section, we prove Theorem \ref{the:Chern} which is a slight improvement
of results in~\cite{Lueck(2002b)}.

In the previous section we have verified that the assumptions of Theorem~0.1 and of
Theorem~0.2 in \cite{Lueck(2002b)} are satisfied in the case where the equivariant homology theory $\calh^?_\ast ( - )$
is given, at a group $G$, by $H_\ast^G ( - ; \JK R ) \otimes_{\IZ} \IQ$, by $H_\ast^G ( - ; \JK^{\minusinfinity} R )
\otimes_{\IZ} \IQ$, or by $H_\ast^G ( - ; \JH\JX^{\otimes_k} R ) \otimes_{\IZ} \IQ$. Let $M$ be a Mackey functor,
for instance $H\mapsto\calh^H_n ( \pt )$ for $n\in\IZ$ fixed. For a finite group $H$, recall the notation
\begin{eqnarray} \label{eq:SC-definition}
\hspace*{3em}
S_H \big( M(H) \big) \; = \;
\coker \left( \bigoplus_{K \lneqq H} \ind_K^H \colon \bigoplus_{K \lneqq H}
M(K) \to
M(H)
\right)
\end{eqnarray}
from \cite{Lueck(2002b)}. Observe for example in the case of $K$-theory that this specializes to \eqref{eq:SC-ist-Artin-defekt}.
We obtain from \cite[Theorems~0.1 and 0.2]{Lueck(2002b)}, for every $G$-$CW$-complex $X$ which is \emph{proper} (i.e.\ with all
stabilizers finite) and every $n\in\IZ$, a canonical isomorphism
\[
\calh_n^G ( X ) \; \cong \; \displaystyle \bigoplus_{p + q = n} \; \displaystyle \bigoplus_{(H) \in (\Fin)} H_p
\big( Z_G H \backslash  X^H ; \IQ \big) \otimes_{\IQ [W_G H ] } S_H \big( \calh^H_{q} ( \pt ) \big)\,,
\]
where $(\Fin)$ denotes the set of conjugacy classes of finite subgroups of $G$. This isomorphism
is natural in $X$ and also in the equivariant homology theory with Mackey structure $\calh^?_\ast ( - )$
(i.e.\ for natural transformations of equivariant homology theories respecting the induction and Mackey
structures). Now, take $X = \underline{E}G$. As in Lemma~8.1 in \cite{Lueck-Reich-Varisco(2003)}, one shows
that the projections
\[
Z_G H \backslash \underline{E}G^H  \leftarrow EZ_G H \times_{Z_G H }  \underline{E}G^H \to EZ_G H / Z_G H = BZ_G H
\]
induce isomorphisms on rational homology. Theorems \ref{the:Chern} now follows
from the next two lemmas.

\begin{lemma} \label{lemma:vanish-not-cyclic}
Let $R$ be a  ring and let $H$ be a finite group. If $H$ is not cyclic, then
\[
S_H \big( K_n(RH) \tensor_{\IZ} \IQ \big) = 0 \quad \mbox{and} \quad S_H \big( \HX^{\otimes_k}_n ( R H )
\otimes_{\IZ} \IQ \big) = 0\,,
\]
for all $n \in \IZ$.
\end{lemma}

\begin{proof}
For a group $H$, let $\Sw (H,\IZ)$ be its Swan group, i.e.\ the Grothendieck group of left $\IZ H$-modules
which are finitely generated as abelian groups. Let $\Sw^f(H,\IZ)$ be the Grothendieck
group of left $\IZ H$-modules which are finitely generated free as abelian groups.
The obvious map $\Sw^f(H,\IZ) \to \Sw(H,\IZ)$ is an isomorphism,
see \cite[Proposition~1.1 on page~553]{Swan(1960a)}.
If $H$ is a finite group, then  $\Sw^f(H ,\IZ)$, and hence also $\Sw(H,\IZ)$, has the structure of
a commutative associative ring, where multiplication is induced by the tensor product over $\IZ$ equipped with the
diagonal $H$-action. The tensor product over $\IZ$ equipped with the diagonal action also leads to a $\Sw^f( H , \IZ)$-module
structure, and hence a $\Sw ( H , \IZ)$-module structure,
on $K_n ( RH )$ for each $n \in \IZ$ and each coefficient ring $R$. For an injective group homomorphism
$\alpha \colon H \into K$ between finite groups, we have the
usual induction and restriction homomorphisms $\ind_H^K \colon \Sw( H , \IZ ) \to \Sw( K , \IZ )$ and
$\res^K_H \colon \Sw( K , \IZ ) \to \Sw( H , \IZ )$. It is not difficult to check that with these structures,
$\Sw (- , \IZ)$ is a Green ring functor with values in abelian groups and that, for each $n \in \IZ$,
the functor $K_n ( R (-) )$ is a module over it (compare~\cite[Sections 7 and 8]{Lueck(2002b)}). Now, by a result
of Swan \cite[Corollary~4.2 on page~560]{Swan(1960a)}, for every finite group $H$, the cokernel of the map
\[
    \bigoplus_{\substack{C\lneqq H \\ C \mbox{ \small cyclic}}}
     \ind_C^H\colon \bigoplus_{\substack{C\lneqq H \\ C \mbox{ \small cyclic}}}
     \Sw(C,\IZ) \tensor_{\IZ}\IQ \to \Sw(H , \IZ) \tensor_{\IZ}\IQ
\]
is annihilated by $|H|^2$. With suitable elements $x_C \in \Sw ( C , \IZ )$, we can hence write
\[
|H|^2 \cdot [ \IZ ] = \sum_{\substack{C\lneqq H \\ C \mbox{ \small cyclic}}} \ind_C^H (x_C )\,.
\]
Therefore, up to multiplication by $|H|^2$, every element $y \in K_n ( R H )$ is induced from proper cyclic subgroups, since
\[
|H|^2 \cdot y = |H|^2 \cdot [ \IZ ] \cdot y =
\sum_{\substack{C\lneqq H \\ C \mbox{ \small cyclic}}} \ind_C^H (x_C ) \cdot y  =
\sum_{\substack{C\lneqq H \\ C \mbox{ \small cyclic}}} \ind_C^H (x_C \cdot \res_C^H y)\,.
\]
The argument for Hochschild homology and its cyclic variants is similar. The module structure over the Swan
ring is also in that case induced by the tensor product over $\IZ$.
\end{proof}

\begin{remark}
More generally, the proof of Lemma~\ref{lemma:vanish-not-cyclic} works for every module over the rationalized
Swan group $\Sw( - , \IZ)\otimes_{\IZ} \IQ$ considered as a Green ring functor. Note that such a statement does
not hold in general for modules over the rationalized Burnside ring $A( - )\otimes_{\IZ} \IQ$ viewed as a Green
ring functor.
\end{remark}

\begin{lemma} \label{lem:module-over-burnside}
Let $C$ be a finite cyclic group and $M$ a Mackey functor with values in $\IQ$-modules.
Keeping notation as in \eqref{theta_C} and \eqref{eq:SC-definition}, there is
a natural isomorphism
\[
\theta_{C} \big( M ( C ) \big) \; \cong \; S_C \big( M ( C ) \big)\,.
\]
\end{lemma}

\begin{proof} Let $D$ be a subgroup of $C$. The ring homomorphism $\chi_C$ of \eqref{def-chi-G}
sends $[C /D]$ to $(x_{(E)})_{(E)}$, where $x_{(E)} = |(C/D)^{E}|$ and hence $x_{(E)}=[C \colon D]$ if $E \leq D$ and
is $0$ otherwise. Therefore the maps $i_D^C$ and $r^C_D$ which make the diagrams
\[
\xymatrix{
A(D) \ar[d]_-{\chi_{D}} \ar[r]^-{\ind_D^C} &  A(C) \ar[d]^-{\chi_{C}}  &
A(C) \ar[d]_-{\chi_{C}} \ar[r]^-{\res^C_D} &  A(D) \ar[d]^-{\chi_{D}} \\
\prod_{\sub D} \IQ  \ar[r]^-{i_D^C} & \prod_{\sub C} \IQ  & \prod_{\sub C} \IQ  \ar[r]^-{r^C_D} &  \prod_{\sub D} \IQ
         }
\]
commutative are easily seen to be given as follows. The map $i_D^C$
is multiplication by the index $[C \colon D]$ followed by the inclusion of the factors corresponding to subgroups of
$C$ contained in $D$. The map $r^C_D$ is the projection onto the factors corresponding to subgroups of $C$ contained
in $D$. In particular, $\chi_C(\ind_D^C (\theta_D))$, considered as a function, is supported only on $(D)$ and takes
there the value $[C \colon D]$. As a consequence, in $A(C) \otimes_{\IZ} \IQ$, we have
\[
1 = [C/C] = \sum_{D \leq C} \frac{1}{[C \colon D]} \ind_D^C ( \theta_D )\,.
\]
Each element in the image of the map $1- \theta_C \colon M(C) \to M(C)$ lies in the image of
$I=\bigoplus_{D \lneqq C} \ind_C^D$, because
\[
(1- \theta_C) x = \left( \sum_{D \lneqq C} \frac{1}{[C \colon D]} \ind_D^C  \theta_D  \right) x
= \sum_{D \lneqq C} \frac{1}{[C \colon D]} \ind_D^C  ( \theta_D  \res^C_D x )\,.
\]
Moreover $\theta_C \colon M(C) \to M(C)$ vanishes on the image of this map $I$; indeed, for $D \lneqq C$,
it follows from the description of $r^C_D$ that $\res^C_D ( \theta_C ) =0$, and therefore
\[
\theta_C  \ind_D^C y  = \ind_D^C ( \res^C_D (\theta_C) y ) = \ind_D^C ( 0 \cdot y ) = 0\,.
\]
So, the cokernel $S_C(M(C))$ of $I$ is isomorphic to the image $\theta_C(M(C))$ of $\theta_C$.
\end{proof}

%%%%%%%%%%%%%%%%%%%%%%%%%%%%%%%%%%%%%%%%%%%%%%%%%%%%%%%%%%%%%%%%%%%%%%%%%%%%%%%%%%%%%%%%%%%%%%%%%%%%%%
%%%%%%%%%%%%%%%%%%%%%%%%%%%%%%%%%%%%%%%%%%%%%%%%%%%%%%%%%%%%%%%%%%%%%%%%%%%%%%%%%%%%%%%%%%%%%%%%%%%%%%
\section{Comparing different models}\label{sec:comparing-different-models}

In order to prove splitting results in Section~\ref{sec:Splitting-assembly-maps},
we will work with a chain complex version and  occasionally with a simplicial abelian
group version of the equivariant homology theory that is associated to Hochschild homology.
In the present section, we define these versions and  prove that they all agree.

%%%%%%%%%%%%%%%%%%%%%%%%%%%%%%%%%

\smallskip

Again, fix a group $G$.
A construction analogous to the balanced smash product \eqref{balanced-smash}, but with smash products ``\,$\wedge$\,''
replaced by tensor products over $\IZ$, and with wedge sums ``\,$\vee$\,'' replaced by direct sums, yields the notion
of \emph{balanced tensor product} $M(?) \otimes_{\IZ \Or G} N(?)$ of a co- and a contravariant $\IZ
\Or G$-module $M(?)$ and $N(?)$. Here by definition a co- or contravariant \emph{$\IZ
\Or G$-module} is a co- respectively  contravariant functor
from $\Or G$ to abelian groups.
Let $C_{\ast} =
C_{\ast} ( ? )$ be a covariant \emph{$\IZ \Or G$-chain complex}, i.e.\ a functor from the orbit category to the category
of chain complexes of abelian groups. We define the \emph{$G$-equivariant Bredon hyperhomology} of a pair of
$G$-$CW$-complexes $(X,A)$ with coefficients in $C_{\ast}$ as
\begin{eqnarray*}
H_\ast^G ( X, A ; C_{\ast}) & = & H_\ast \Big( \Tot^{\oplus} \big( \widetilde{C}_{\ast}^{\sing}( (X_+/A_+)^? )
\otimes_{\IZ \Or G} C_{\ast}(?) \big) \Big)\,.
\end{eqnarray*}
Here, for a pointed $G$-space $Y=(Y,y_0)$ (where $y_0$ is a
$G$-fixpoint), the functor which sends $G/H$ to the reduced
singular chain complex of $Y^H$ is denoted by
$\widetilde{C}_{\ast}^{\sing}(Y^?)$. In this construction, up to
canonical isomorphism, we can replace
$\widetilde{C}_{\ast}^{\sing}(-)$ by the \emph{reduced cellular
chain complex} $\widetilde{C}_{\ast}^{\cell}(-)$ (this will be
needed in Subsection~\ref{subsec:splitting-HN-HP-assembly}).

For a \emph{simplicial $\IZ\Or G$-module} $M_{\bullet}=M_{\bullet} (?)$, i.e.\ a covariant functor from $\Or G$
to the category of simplicial abelian groups, we define similarly
\begin{eqnarray*}
H_{\ast}^G ( X, A ; M_{\bullet}) & = & \pi_{\ast} \Big( \big|\widetilde{\IZ} \big[ S_{\bullet} ( (X_+/A_+)^? ) \big]
\otimes_{\IZ \Or G} M_{\bullet} ( ? ) \big| \Big)\,.
\end{eqnarray*}
Here, $S_{\bullet}$ stands for the singular simplicial set associated to a topological space. For a pointed
simplicial set $Y_\bullet=(Y_\bullet, y_0)$, we set $\widetilde{\IZ} [Y_\bullet] = \IZ [Y_\bullet] \big/ \IZ [y_0]$
and the tensor
products of simplicial abelian groups are taken degreewise.

For an $\Or G$-spectrum $\JE=\JE(?)$, recall that we use the notation
\begin{eqnarray*}
H_{\ast}^G( X,A ; \JE ) & = & \pi_{\ast} \big( (X_+/A_+)^? \sma_{\Or G} \JE(?) \big)\,.
\end{eqnarray*}

A simplicial $\IZ\Or G$-module $M_{\bullet}= M_{\bullet}(?)$ gives rise to a $\IZ \Or G$-chain complex
$\DK_{\ast}\!M_{\bullet}$ {\sl via} the Dold-Kan correspondence and determines an $\Or G$-spectrum
$\JH M_{\bullet}$ {\sl via} the Eilenberg-Mc\,Lane functor (see Subsections~\ref{subsec:simpl-dold-kan}
and \ref{subsec:Gamma-spaces}). The following proposition specializes to a well-known fact in the
case where $G$ is the trivial group.

\begin{proposition} \label{prop:compare}
Let $M_{\bullet}$ be a functor from $\Or G$ to simplicial abelian groups. Then, there are natural
isomorphisms of $G$-homology theories defined on pairs of $G$-$CW$-complexes
\begin{eqnarray*}
H_\ast^G ( - ; \JH M_{\bullet} )  \xrightarrow{\cong}
H_\ast^G( - ; M_{\bullet} )
\xrightarrow{\cong}
H_\ast^G( - ; \DK_{\ast}\!M_{\bullet} )\,.
\end{eqnarray*}
\end{proposition}

In particular we have natural isomorphisms
\begin{eqnarray} \label{some-isos}
H_\ast^G ( - ; \JH \JX^{\otimes_k}R ) \; \cong \; H_\ast^G( - ; C_{\bullet}^{\HX^{\otimes_{\!k}}\!R}) \; \cong \;
H_\ast^G( - ; C_{\ast}^{\HX^{\otimes_{\!k}}\!R} )\,.
\end{eqnarray}
Here we have used the notation
\begin{eqnarray}
C_{\bullet}^{\HX^{\otimes_{\!k}}\!R}(? ) & = &  C_{\bullet}^{\HX^{\otimes_{\!k}}} R \calg^G(?) \label{def-C-bullet-HXkR} \\
C_{\ast}^{\HX^{\otimes_{\!k}}\!R}(? ) & = & C_{\ast}^{\HX^{\otimes_{\!k}}} R \calg^G(? )\,, \label{def-C-ast-HXkR}
\end{eqnarray}
for the indicated simplicial $\IZ \Or G$-module respectively $\IZ \Or G$-chain complex,
compare Subsection~\ref{define-bicomplexes}, and \eqref{def-HX-OrG-spectrum}.

\begin{proof}[Proof of Proposition~\ref{prop:compare}]
We discuss the first natural transformation in the absolute case, i.e.\ for $A= \varnothing$ (the general
case is similar).
For a spectrum $\JE$ in simplicial sets we denote by $| \JE |$ the associated (topological) spectrum.
For every $G$-$CW$-complex $X$ and every $\Or G$-spectrum $\JE=\JE(?)$ in simplicial sets,
there is a natural equivalence and, since realization commutes with taking coequalizers, a natural homeomorphism
\[
X_+ \sma_{\Or G} | \JE |
\xleftarrow{\simeq}
|S_{\bullet} X_+ | \sma_{\Or G} | \JE |
\xrightarrow{\cong}
\big|S_{\bullet} X_+ \sma_{\Or G} \JE \big|\,.
 \]
Note that there is an obvious natural isomorphism of spectra
$$
S_{\bullet} X_+ \sma_{\Or G}\IH M_{\bullet}(\JS) \; \cong \; \Big(S_{\bullet} X_+ \sma_{\Or G}\IH M_{\bullet}\Big)(\JS)\,.
$$
Observe also that for an (unpointed) simplicial set $Y_\bullet$, there is an isomorphism of $\Gamma$-spaces
\begin{eqnarray} \label{easy-iso-for-prop-compare}
\IZ[Y_\bullet] \otimes_{\IZ \Or G} \IH M_{\bullet} \; \cong \;
\IH \big( \IZ[Y_\bullet] \otimes_{\IZ \Or G} M_{\bullet} \big)\,.
\end{eqnarray}
By \eqref{homotopy-HM}, the homotopy groups of the spectrum associated to the right hand side are given by the (unstable)
homotopy groups of (the realization of) $\IZ[Y_\bullet] \otimes_{\Or G} M_{\bullet}$\,.

So, observing that $\widetilde{\IZ}
(S_\bullet X_+) \cong \IZ[S_\bullet X]$, to produce the first natural transformation of the statement, it will suffice
to define a natural transformation of $\Gamma$-spaces
\begin{eqnarray} \label{nat-trans-for-prop-compare}
S_{\bullet} X_+ \sma_{\Or G} \IH M_{\bullet} \to \widetilde{\IZ}[ S_{\bullet} X_+ ] \otimes_{\IZ \Or G} \IH M_{\bullet}\,.
\end{eqnarray}
More generally, for every contravariant functor $Z_\bullet=Z_\bullet(?)$ from $\Or G$ to pointed simplicial sets and
every covariant functor $N_{\bullet}=N_{\bullet}(?)$ from $\Or G$ to simplicial abelian groups,
we will construct a natural transformation
\[
Z_{\bullet} \sma_{\Or G}  N_{\bullet} \to
\widetilde{\IZ} [ Z_{\bullet} ] \otimes_{\IZ \Or G}  N_{\bullet}\,.
\]
To produce the map we use the following facts. The left-hand side is defined as a coequalizer in the category of pointed
simplicial sets completely analogous to \eqref{balanced-smash}, and the right-hand side similarly as a coequalizer in the category of
simplicial abelian groups.
Let $U$ denote the forgetful functor from simplicial abelian groups to pointed simplicial sets.
For a pointed simplicial set $X$, a simplicial abelian group $A$ and a family $A_i$, $i \in I$ of simplicial abelian groups there are
obvious natural maps $X \sma UA \to U ( \widetilde{\IZ} [ X ] \otimes A )$ and $\bigvee_{i \in I} UA_i \to U(\bigoplus_{i \in I} A_i)$.
Given two maps $f ,g \colon A \to B$ of simplicial abelian groups
there is an obvious natural map $\coequ( Uf , Ug ) \to U \coequ(f,g)$.
Combining these facts one easily constructs the required natural transformation above.

Now, we show that the first natural transformation of the statement is an isomorphism. If a natural transformation between
$G$-homology theories induces an isomorphism when evaluated on all orbits $G/H$, then it induces an isomorphism for all pairs
of $G$-$CW$-complexes by a well-known argument. Unravelling the construction of the first natural transformation, it suffices
to check that for every orbit $G/H$, the map
$$
%S_{\bullet} (G/H)_+ \sma_{\Or G}\IH M_{\bullet}(\JS) =
\Big(S_{\bullet} (G/H)_+ \sma_{\Or G}\IH M_{\bullet}\Big)(\JS) \to
\IH \big( \widetilde{\IZ}[S_{\bullet} (G/H)_+] \otimes_{\IZ \Or G} M_{\bullet} \big)(\JS)
$$
induced by \eqref{easy-iso-for-prop-compare} and \eqref{nat-trans-for-prop-compare} is a stable weak equivalence.
Before evaluation at $\JS$, both sides are canonically isomorphic to the $\Gamma$-space $\IH M_{\bullet}(G/H)$
by suitable analogues of Lemma~\ref{coyoneda}.
We leave it to the reader to verify that we indeed have $G$-homology
theories here, compare \cite[Lemma~4.2]{Davis-Lueck(1998)}.

We now construct the second natural transformation of the
statement and prove at the same time that it is an isomorphism.
For a bisimplicial abelian group $A_{\bullet \bullet}$, let
$C_{\ast\ast} ( A_{\bullet \bullet} )$ denote the associated
bicomplex, compare \cite[page~275]{Weibel(1994)}. Note that given
two simplicial abelian groups $C_\bullet$ and $D_\bullet$, there
is a natural isomorphism of bicomplexes $\DK_\ast (C_\bullet)
\otimes_\IZ \DK_\ast(D_\bullet) \cong C_{\ast \ast} (C_\bullet
\otimes_\IZ D_\bullet)$, where $C_\bullet \otimes_\IZ D_\bullet$
is viewed as a bisimplicial abelian group. Note also that
$\DK_\ast \big( \widetilde{\IZ}[ S_\bullet (X) ]\big) =
\widetilde{C}^{\sing}_\ast (X)$, for every space $X$. The
degreewise tensor products of simplicial abelian groups appearing
in the source of the second natural transformation can be thought
of as diagonals of the corresponding bisimplicial sets. Applying
all these observations and using again the definition of the
balanced tensor product in terms of coequalizers, it suffices to
observe that for every pair of maps $f_{\bullet \bullet},
g_{\bullet \bullet} \colon A_{\bullet \bullet} \to B_{\bullet
\bullet}$ of bisimplicial abelian groups, we have the following
chain of isomorphisms
\begin{eqnarray*}
\pi_\ast \big( \big| \coequ ( \diag f_{\bullet \bullet} , \diag g_{\bullet \bullet}) \big| \big)
& \cong & \pi_\ast \big( \big| \diag \coequ ( f_{\bullet \bullet} , g_{\bullet \bullet}) \big| \big) \\
& \cong & H_\ast \big( \Tot^{\oplus} C_{\ast \ast} (\coequ ( f_{\bullet \bullet} , g_{\bullet \bullet}) ) \big) \\
& \cong & H_\ast \big( \coequ ( \Tot^{\oplus} C_{\ast \ast} ( f_{\bullet \bullet}) , \Tot^{\oplus} C_{\ast \ast}( g_{\bullet \bullet}) ) \big)\,,
\end{eqnarray*}
where the second isomorphism is the Eilenberg-Zilber Theorem as formulated in \cite[Theorem~8.5.1 on page~276]{Weibel(1994)}.
\end{proof}

%%%%%%%%%%%%%%%%%%%%%%%%%%%%%%%%%%%%%%%%%%%%%%%%%%%%%%%%%%%%%%%%%%%%%%%%%%%%%%%%%%%%%%%
%%%%%%%%%%%%%%%%%%%%%%%%%%%%%%%%%%%%%%%%%%%%%%%%%%%%%%%%%%%%%%%%%%%%%%%%%%%%%%%%%%%%%%%
\section{Splitting assembly maps}\label{sec:Splitting-assembly-maps}

In this section, we prove Theorem~\ref{the:splitHochschild} and Addendum~\ref{add:split-HN},
i.e.\ the splitting and isomorphism results for the assembly maps in Hochschild, cyclic,
periodic cyclic and negative cyclic homology. We begin with the case of Hochschild homology.

%%%%%%%%%%%%%%%%%%%%%%%%%%%%%%%%%%%%%%%%%%%%%%%%%%%%%%%%%%%%%%%%%%%%%%%%%%%%%%%%%%%%%%%
\subsection{Splitting the Hochschild homology assembly map} \label{subsec:splitting-HH}

Fix a group $G$ and let $S$ be a $G$-set. Recall that $\con G$ denotes the set of conjugacy classes of $G$.
Sending  a $q$-simplex
$(g_0, \ldots , g_q)$ in $\CN_{\bullet} \calg^G(S)$ to the conjugacyclass $\conjclass{g_0 \cdots g_q}$
yields a map of cyclic sets
\[
\CN_{\bullet} \calg^G(S) \to \con G\,.
\]
Here $\con G$ is considered as a  constant cyclic set.
The cyclic nerve decomposes, as a cyclic set, into the disjoint union of the corresponding pre-images, namely
\begin{eqnarray}\label{decomposition}
\CN_{\bullet} \calg^G (S) \; = \; \coprod_{\conjclass{c} \in \con G } \CN_{\bullet \conjclass{c} } \calg^G (S)\,.
\end{eqnarray}
Observe that $\CN_{\bullet \conjclass{c} } \calg^G (G/H) \neq \varnothing$ implies that $\langle c \rangle$ is
subconjugate to $H$.
For every small category $\calc$, we have a natural isomorphism
\[
k\CN_{\bullet} \calc \; \cong \; \CN_{\bullet}^{\otimes_k}   k \calc
\]
and, because of the isomorphism $R \calc \cong R \otimes_k k \calc$, also
\begin{eqnarray} \label{decompose-k-R}
\CN_{\bullet}^{\otimes_k} R \calc \; \cong \; (\CN_{\bullet}^{\otimes_k} k \calc ) \otimes_k ( \CN_{\bullet}^{\otimes_k} R)\,.
\end{eqnarray}
We therefore obtain an induced  decomposition for the $k$-linear cyclic nerve of $R \calg^G (S)$, that we denote by
\begin{eqnarray} \label{linear-splitting}
\CN_{\bullet}^{\otimes_k} R \calg^G (S) & = &
\bigoplus_{\conjclass{c} \in \con G } \CN_{\bullet \conjclass{c}}^{\otimes_k} R \calg^G (S) \\
& = &
\bigoplus_{\conjclass{c} \in \con G } \CN_{\bullet \conjclass{c} }^{\otimes_k}  k\calg^G (S) \otimes_k \CN_{\bullet}^{\otimes} R\,.
\label{decompositionlinear}
\end{eqnarray}
For typographical reasons we introduce the following abbreviation
for the corresponding decomposition of simplicial $\IZ\Or G$-modules\,:
\begin{eqnarray} \label{decomposition of C_bullet^HH}
C_{\bullet}^{\HH^{\otimes_{\!k}}\!R} (?) & = &  \bigoplus_{\conjclass{c} \in \con G }
C_{\bullet \conjclass{c}}^{\HH^{\otimes_{\!k}}\!R} (?)\,;
\end{eqnarray}
see \eqref{def-C-bullet-HXkR} for the notation. Using the identifications \eqref{some-isos}
and the decomposition \eqref{decomposition of C_bullet^HH}, the Hochschild homology generalized assembly map
$$
\xymatrix@C=4em{
H_n^G(\EGF{G}{\calf};\JH \JH^{\otimes_k}R) \ar[r]^-{\assembly} &
H_n^G(\pt;\JH \JH^{\otimes_k}R) \cong \HH_n^{\otimes_k}(RG)
               }
$$
appearing in diagram \eqref{diagram-assembly-dtr} can be identified with the upper horizontal map in the following
commutative diagram
\begin{eqnarray} \label{proof-diagram}
\hspace*{2.3em}
\vcenter{\xymatrix@C=3.8em{
H_{n}^G\big(\EGF{G}{\calf} ; \bigoplus_{\conjclass{c} \in \con G} C_{\bullet \conjclass{c}}^{\HH^{\otimes_{\!k}}\!R}\big)
\ar[d]_{\id \otimes pr_{\calf}} \ar[r]^-{\assembly}
&
H_{n}^G\big( \pt  ; \bigoplus_{\conjclass{c} \in \con G} C_{\bullet \conjclass{c}}^{\HH^{\otimes_{\!k}}\!R}\big)
\ar[d]^{\id \otimes pr_{\calf}}
\\
H_{n}^G\big(\EGF{G}{\calf} ; \smash[b]{\bigoplus_{\stackrel{\conjclass{c} \in \con G}{\langle c \rangle \in \calf}}}
C_{\bullet \conjclass{c}}^{\HH^{\otimes_{\!k}}\!R}\big) \ar[r]
&
H_{n}^G\big( \pt ; \smash[b]{\bigoplus_{\stackrel{\conjclass{c} \in \con G}{\langle c \rangle \in \calf}}}
C_{\bullet \conjclass{c}}^{\HH^{\otimes_{\!k}}\!R}\big).
         }}
\end{eqnarray}

\smallskip

\noindent
Here, the vertical maps are induced by the projection $pr_{\calf}$ onto the summands for which the cyclic subgroup
$\langle c \rangle$ belongs to the family $\calf$. Note that $pr_{\calf}$ is the identity map if $\calf$ contains
all cyclic subgroups of $G$. The statement about Hochschild homology in Theorem~\ref{the:splitHochschild} now follows
directly from the following two lemmas.

\begin{lemma} \label{leftvertical}
For every family $\calf$, the left vertical map in \eqref{proof-diagram} is an isomorphism.
\end{lemma}

\begin{lemma} \label{lowerhorizontal}
For every family $\calf$, the bottom map in \eqref{proof-diagram} is an isomorphism.
\end{lemma}

The proofs of Lemmas~\ref{leftvertical} and \ref{lowerhorizontal} will occupy the rest of this subsection. They rely on
the following computation of the cyclic nerve of a transport groupoid.

\smallskip

Let $E_{\bullet} G$ be the simplicial set given by $N_{\bullet} \calg^G( G/1 )$. In words\,: consider $G$ as a category with
$G$ as set of objects and precisely one morphism between any two objects, and then take the nerve of this category.
This is a simplicial model for the universal free $G$-space which is usually denoted by $EG$.
For $c \in G$ let $\langle c \rangle$ be the cyclic subgroup generated by $c$.
For $h \in N_G \langle c \rangle$, let $R_h \in \map_G( G / \langle c \rangle, G/ \langle c \rangle )$
be the map given by $R_h ( g \langle c \rangle ) = gh \langle c \rangle$. For every $G$-set $S$,
precomposing with
$R_h$ yields a left action of $Z_G \langle c \rangle \leq N_G \langle c \rangle$
on $\map_G (G/\langle c \rangle , S)$.

\begin{proposition} \label{keyproposition}
For a group $G$, choose a representative $c \in \conjclass{c}$ for each conjugacy class $\conjclass{c} \in \con G$.
Let $\langle c \rangle$ denote the cyclic subgroup it generates.
There is a map of $\Or G$-simplicial sets (depending on the choice)
\[
\coprod_{\conjclass{c} \in \con G}
E_{\bullet} Z_G \langle c \rangle
\times_{Z_G \langle c \rangle}
\map_{G}( G/ \langle c \rangle , ? )
\to \CN_{\bullet} \calg^{G}( ? )\,.
\]
This map is objectwise a simplicial homotopy equivalence, and is compatible with the decomposition
\eqref{decomposition} of the target.
\end{proposition}

\begin{remark}
There seems to be no obvious cyclic structure on the source of the map above.
\end{remark}

\begin{proof}[Proof of Proposition \ref{keyproposition}]
We first introduce some more notation. Given a groupoid $\calg$, we denote by $\aut \calg$ its
category of automorphisms, i.e.\ the category whose objects are automorphisms $h\colon s \to s$
in $\calg$ and where a morphism from $h\colon s \to s$ to $h'\colon t \to t$
is given by a morphism $g\colon s \to t$ satisfying $h' \circ g = g \circ h$. In the case where
$\calg=\calg^G (S)$, the conjugacy class $\conjclass{h} \in \con G$ associated to an object $h\colon s
\to s$ in $\aut \calg^G ( S )$ does only depend on the isomorphism class of this object.
This yields a well-defined map of simplicial sets
\[
N_{\bullet} \aut \calg^G (S) \; \longrightarrow \; \con G\;,\qquad
\vcenter{\xymatrix@R=1.5em@C=1.9em{
s_0 \ar[d]_{h_0} & s_1  \ar[d]_{h_1} \ar[l]_-{g_0} & \downdots \ar[l]_-{g_1} & s_q \ar[d]_{h_q} \ar[l]_-{\,g_{q-1}} \\
s_0 & s_1 \ar[l]^-{g_0} & \downdots \ar[l]^-{g_1} & s_q \ar[l]^-{\,g_{q-1}}.
         }}
\quad\longmapsto\quad\conjclass{h_0}\,,
\]
where $\con G$ is considered as a constant simplicial set.
Let $\aut_{\conjclass{c}} \calg^G (S)$ denote the full subcategory of $\aut \calg^G (S)$ on the objects $h\colon s
\to  s$ with $h \in \conjclass{c}$.
The decomposition of the nerve into pre-images
under the map to $\con G$ above is given by
\[
N_{\bullet} \aut \calg^G (S) \; = \; \coprod_{\conjclass{c} \in \con G }  N_{\bullet} \aut_{\conjclass{c}} \calg^G (S )\,.
\]
The components of the map in Proposition~\ref{keyproposition} are obtained
as the
composition of the three maps
\[
\xymatrix{
E_{\bullet} Z_G \langle c \rangle
\times_{Z_G \langle c \rangle}
\map_G( G/ \langle c \rangle , ? )
\ar[r] &
N_{\bullet} \calg^{Z_G \langle c \rangle} \big( \map_G ( G/ \langle c \rangle , ? ) \big)  \ar[dl] \\
N_{\bullet} \aut_{\conjclass{c}} \calg^G ( ? )
\ar[r] & \CN_{\bullet \conjclass{c}} \calg^G( ? )
         }
\]
which are constructed in
the following lemma.  Proposition~\ref{keyproposition} is an immediate consequence of that lemma.
\end{proof}

\begin{lemma}
Let $G$ be a group and $S$ a $G$-set.
\begin{enumerate}
\item
There is a simplicial isomorphism $E_{\bullet} G \times_G S \to N_{\bullet} \calg^G (S)$.
\item
For $\conjclass{c} \in \con G$, choose a representative $c \in \conjclass{c}$. There is an equivalence of categories
\[
\calg^{Z_G \langle c \rangle} \big( \map_G ( G / \langle c \rangle , S ) \big) \to \aut_{\conjclass{c}}
\calg^G ( S )\,,
\]
which depends on the choice.
\item
For every groupoid $\calg$, there is a simplicial isomorphism
\[
N_{\bullet} \aut \calg  \to \CN_{\bullet} \calg\,.
\]
If $\calg=\calg^G (S)$ then the isomorphism commutes with the maps to $\con G$.
\end{enumerate}
All three constructions are natural with respect to $S$.
\end{lemma}

\begin{proof}
(i) The isomorphism $E_\bullet G \times_G S \to N_\bullet \calg^G(S)$ is given, on level $q$, by
\[
\xymatrix@C=2.5em{
\Big[\, g_0  &  g_1 \ar[l]_-{\;g_0 g_1^{-1} } & \downdots \ar[l]_-{\;g_1 g_2^{-1}} & g_q \ar[l]_-{\;g_{q-1} g_q^{-1}}
\, , \, s \,\Big]
\;\;\longmapsto\;\;
\Big(g_0 s &  g_1 s \ar[l]_-{\;g_0 g_1^{-1} } & \downdots \ar[l]_-{\;g_1 g_2^{-1}} & g_q s \ar[l]_-{\;g_{q-1} g_q^{-1}}\Big)\,.
         }
\]

\smallskip

(ii) The functor sends an object $\phi \in \map_G ( G/\langle c \rangle , S )$ to the automorphism
$c \colon  \phi(e \langle c \rangle ) \to \phi(e \langle c \rangle ) $. Here $e$ is the trivial
element in $G$ and $c$ the chosen representative in $\conjclass{c}$. A morphism $z \colon \phi \to z \phi$, with
$z\in Z_G \langle c \rangle$, is taken to the (iso)morphism $z^{-1} \colon \phi(e \langle c \rangle ) \to z^{-1} \phi
( e \langle c \rangle )$. The functor is full and faithful and every object
in the target category is isomorphic to an image object.

\smallskip

(iii) The isomorphism $N_\bullet \aut \calg  \to \CN_\bullet \calg$ is given, on level $q$, by
\[
\;\;\vcenter{\xymatrix@R=1.5em@C=1.9em{
s_0 \ar[d]_{h_0} &
s_1  \ar[d]_{h_1} \ar[l]_-{g_0} &
\downdots  \ar[l]_-{g_1} &
s_q \ar[d]_{h_q} \ar[l]_-{\,g_{q-1}} &
\ar@{}[d]^{\longmapsto} \\
s_0 &
s_1 \ar[l]^-{g_0} &
\downdots \ar[l]^-{g_1} &
s_q \ar[l]^-{\,g_{q-1}}  &
         }}
\qquad\;
\vcenter{\xymatrix@R=1.5em@C=1.9em{
\!\!s_0  \ar@(dl,dr)[rrr]_{ h_0(g_0 \cdots g_{q-1})^{-1} }  & s_1 \ar[l]_-{g_0} & \downdots \ar[l]_-{g_1} &
s_q \ar[l]_-{\,g_{q-1}}\,.\!\!\!\!\! \\
         }}
\]
The compatibility with the maps to $\con G$ is clear.
\end{proof}

The following is the linear analogue of Proposition~\ref{keyproposition}.

\begin{corollary} \label{keypropositionlinear}
For every conjugacy class $\conjclass{c} \in \con G$ there is
natural transformation of functors from the orbit category
$\Or G$ to the category of simplicial $k$-modules,
\[
k [E_{\bullet} Z_G \langle c \rangle]
\otimes_{k Z_G\langle c \rangle }
k \map( G/\langle c \rangle , ? )
\otimes_k
\CN_{\bullet}^{\otimes_k} R
\to
\CN_{\bullet \conjclass{c}}^{\otimes_k} R \calg^G( ? )
\]
which is objectwise a homotopy equivalence.
\end{corollary}
\begin{proof}
Apply the functor free $k$-module $k(-)$ to the map in Proposition~\ref{keyproposition} and recall the identification
\eqref{decompose-k-R}.
\end{proof}

Observe that we have a decomposition of $G$-homology theories
\begin{eqnarray} \label{decomp}
H_\ast^G \big( - ; C_{\bullet}^{\HH^{\otimes_{\!k}}\!R} \big)
%\cong H_\ast^G \Big( - ; \bigoplus_{\conjclass{c} \in \con G } C_{\bullet  \conjclass{c}}^{\HH^{\otimes_{\!k}}\!R} \Big)
\cong \bigoplus_{\conjclass{c} \in \con G } H_\ast^G \big( - ; C_{\bullet  \conjclass{c}}^{\HH^{\otimes_{\!k}}\!R} \big)\,,
\end{eqnarray}
because the tensor product over the orbit category and homology both commute with direct sums.
For each of the summands, we have the following computation.

\begin{proposition} \label{computation-one-summand}
For every $G$-$CW$-complex $X$ and every $\conjclass{c} \in \con G$, there is a natural isomorphism
\[
H_{\ast}^G \big( X ; C_{\bullet \conjclass{c}}^{\HH^{\otimes_{\!k}}\!R} \big) \; \cong \;
H_\ast \big( X^{\langle c \rangle} \times_{Z_G \langle c \rangle} E Z_G {\langle c \rangle} ;
\CN_{\bullet}^{\otimes_k} R \big)\,.
\]
\end{proposition}

\begin{proof}
On the level of simplicial abelian groups, Corollary~\ref{keypropositionlinear}, in combination with
Lemma~\ref{coyoneda}, yields
$$
\renewcommand{\arraystretch}{1.7}
\begin{array}{l}
\widetilde{\IZ} [S_{\bullet} X^?_+ ] \otimes_{\IZ \Or G} C_{\bullet \conjclass{c}}^{\HH^{\otimes_{\!k}}\!R} (?) \; \simeq
\qquad\qquad\quad \\
\qquad\qquad\quad \simeq \; \widetilde{\IZ} [ S_{\bullet} X^?_+ ] \otimes_{\IZ \Or G } k[ E_{\bullet} Z_G \langle c \rangle ]
\otimes_{k Z_G \langle c \rangle} k \map ( G / \langle c \rangle , ? ) \otimes_k \CN_{\bullet}^{\otimes_k} R \\
\qquad\qquad\quad \cong \; \widetilde{\IZ} [ S_{\bullet} X^{\langle c \rangle}_+ ] \otimes_{\IZ Z_G \langle c \rangle }
\IZ [ E_{\bullet} Z_G \langle c \rangle ] \otimes_{\IZ} \CN_{\bullet}^{\otimes_k} R\,,
\end{array}
$$
hence the result.
\end{proof}

\begin{lemma} \label{coyoneda}
Let $F$ be a contravariant functor from $\Or G$ to simplicial $k$-modules. Then, for every subgroup $H\leq G$,
there is a natural isomorphism
\[
F(?) \otimes_{k \Or G} k \map_G( G/H , ? ) \; \cong \; F(G/H)\,.
\]
\end{lemma}

We can now finish the proof of the part of Theorem~\ref{the:splitHochschild} concerned with Hochschild homology.

\begin{proof}[Proof of Lemmas~\ref{leftvertical} and \ref{lowerhorizontal}]
Compute the relevant maps in diagram \eqref{proof-diagram} using \eqref{decomp} and Proposition~\ref{computation-one-summand}.
Observe that by the very definition of $\EGF{G}{\calf}$, we have $\EGF{G}{\calf}^{\langle c \rangle} = \varnothing$
if and only if $\langle c \rangle \in \calf$. So the projection $\id \otimes pr_{\calf}$ is the zero map
exactly on those summand which are anyway trivial. This proves Lemma~\ref{leftvertical}. For
$\langle c \rangle \in \calf$, the map $\EGF{G}{\calf}^{\langle c \rangle} \to \pt$ is a homotopy equivalence. Therefore,
$\EGF{G}{\calf}^{\langle c \rangle} \times EZ_G \langle c \rangle \to\pt \times EZ_G \langle c \rangle$ is an equivalence
of free $Z_G \langle c \rangle$-spaces and hence remains an equivalence if we quotient out the $Z_G \langle c \rangle$-action.
This establishes Lemma~\ref{lowerhorizontal}.
\end{proof}

The following example gives a further illustration of the computation achieved above.

\begin{example} \label{rem:comparison_with_standard_splitting}
Combining \eqref{some-isos}, the isomorphism \eqref{decomp} and Proposition~\ref{computation-one-summand},
we get, for every $G$-$CW$-complex $X$, a decomposition
$$H_\ast^G(X;\JH \JH^{\otimes_k}R) \; \cong \;  \bigoplus_{\conjclass{c} \in \con G }
H_\ast\big( X^{\langle c \rangle} \times_{Z_G \langle c \rangle} E Z_G {\langle c \rangle} ;
C_\ast^{\HH^{\otimes_{\!k}}}( R) \big)\,,$$
where each direct summand on the right-hand side is the (non-equivariant) hyper\-homology of the space
$X^{\langle c \rangle} \times_{Z_G \langle c \rangle} E Z_G {\langle c \rangle}$
with coefficients in the Hochschild complex, i.e.\ in the $k$-chain complex
$C_\ast^{\HH^{\otimes_{\!k}}} (R) = \DK_{\ast} ( \CN_{\bullet}^{\otimes_k} R )$.  In the case $R = k$, the degree zero inclusion
$k \to C^{\HH^{\otimes_{\!k}}} (k)$ is a homology equivalence and hence a chain homotopy equivalence,
because both complexes are bounded below and consist of projective $k$-modules. Thus, we infer
\[
H_\ast\big( X^{\langle c \rangle} \times_{Z_G \langle c \rangle} E Z_G {\langle c \rangle} ; C_\ast^{\HH^{\otimes_{\!k}}} (k) \big)
\; \cong \; H_\ast\big(X^{\langle c \rangle} \times_{Z_G \langle c \rangle} E Z_G {\langle c \rangle};k\big)
\]
for each conjugacy class $\conjclass{c}$, and therefore
$$H_\ast^G \big( X;\JH \JH^{\otimes_k}R \big) \; \cong \;
 \bigoplus_{\conjclass{c} \in \con G }
H_\ast \big( X^{\langle c \rangle} \times_{Z_G \langle c \rangle} E Z_G {\langle c \rangle};k \big)\,.$$
In the special case where $X = \pt$ and $R = k$, we rediscover the well-known decomposition
of $k$-modules
\begin{eqnarray}
\label{standard-decomposition_of_HH(KG)}
\HH^{\otimes_k}(kG) &  \cong &
 \bigoplus_{\conjclass{c} \in \con G }  H_\ast \big( B Z_G {\langle c \rangle};k \big)\,.
\end{eqnarray}
If we insert $X = \EGF{G}{\calf}$ for an arbitrary family of subgroups $\calf$, we obtain
\[
H_\ast^G \big( \EGF{G}{\calf};\JH \JH^{\otimes_k}R \big) \; \cong \;
 \bigoplus_{ \substack{\conjclass{c} \in \con G \\ \langle c \rangle \in \calf}}  H_\ast \big( B Z_G {\langle c \rangle};k \big)\,,
\]
because $\EGF{G}{\calf}^{\langle c \rangle} \times EZ_G \langle c \rangle$ is a model for $E Z_G \langle c \rangle$ if
$\langle c \rangle \in \calf$ and is empty otherwise. The map $\EGF{G}{\calf} \to \pt$ induces the obvious inclusion.
\end{example}

\begin{remark}
Of course one does not need the elaborate setup using spectra, nor Theorem~\ref{the:splitHochschild},
in order to prove the well-known decomposition \eqref{standard-decomposition_of_HH(KG)}. But our aim was to compare
the Hochschild assembly map with the one for K-theory. There is no chain complex version of the assembly map
on the level of $K$-theory. Furthermore, our effort has the pay off that it can be generalized to topological
Hochschild homology and its refinements as explained in \cite{Lueck-Reich-Rognes-Varisco(2005)}.
\end{remark}

%%%%%%%%%%%%%%%%%%%%%%%%%%%%%%%%%%%%%%%%%%%%%%%%%%%%%%%%%%%%%%%%%%%%%%%%%%%%%%%%%%%%%%%
\subsection{Splitting cyclic, periodic cyclic and negative cyclic assembly maps} \label{subsec:splitting-HN-HP-assembly}

Observe that the sum decomposition \eqref{linear-splitting} is compatible with the cyclic structure.
Keeping notation as in \eqref{def-C-ast-HXkR}, we hence obtain a decomposition
\begin{eqnarray*}
C_{\ast}^{\HX^{\otimes_{\!k}}\!R}(?) & = & \bigoplus_{\conjclass{c} \in \con G }
C_{\ast \conjclass{c}}^{\HX^{\otimes_{\!k}}\!R}(?)
\end{eqnarray*}
of $\IZ\Or G$-chain complexes. Compare with the splitting \eqref{decomposition of C_bullet^HH}.
There is consequently a version of diagram \eqref{proof-diagram}
with $C_{\bullet \conjclass{c}}^{\HH^{\otimes_{\!k}}\!R}$ replaced everywhere by $C_{\ast \conjclass{c}}^
{\HX^{\otimes_{\!k}}\!R}$, and where the upper horizontal map corresponds to the generalized assembly map
for $\HX$-homology.
In order to prove the cyclic homology
part of Theorem~\ref{the:splitHochschild} and to establish Addendum~\ref{add:split-HN},
it suffices to obtain the analogues of Lemmas~\ref{leftvertical} and \ref{lowerhorizontal} with $C_{\bullet
\conjclass{c}}^{\HH^{\otimes_{\!k}}\!R}$ replaced everywhere by $C_{\ast \conjclass{c}}^{\HX^{\otimes_{\!k}}\!R}$.
However, this follows immediately from the following proposition.

\begin{proposition}
Let $X \to X^{\prime}$ be a map of $G$-$CW$-complexes and $Z_{\bullet}(?) \to Z_{\bullet}^{\prime}(?)$
a map of \emph{cyclic $\IZ\Or G$-modules}, i.e.\ a natural transformation between functors from the orbit
category $\Or G$ to the category of cyclic abelian groups. Keep notation as at the
beginning of Subsection~\ref{define-bicomplexes}.

If the induced map
\begin{eqnarray*}
H^G_{\ast} \big( X ; C_{\ast}^{\HH} ( Z_{\bullet} ) \big) \to H^G_{\ast} \big( X^{\prime} ;
C_{\ast}^{\HH} ( Z_{\bullet}^{\prime} ) \big)
\end{eqnarray*}
is an isomorphism, then the map
\begin{eqnarray*}
H^G_{\ast} \big( X ; C_{\ast}^{\HC} ( Z_{\bullet} ) \big) \to
H^G_{\ast} \big( X^{\prime} ; C_{\ast}^{\HC} ( Z_{\bullet}^{\prime} ) \big)
\end{eqnarray*}
is an isomorphism. If moreover $X$ and $X^{\prime}$ are finite $G$-$CW$-complexes, then also
the maps
\begin{eqnarray*}
H^G_{\ast} \big( X ; C_{\ast}^{\HP} ( Z_{\bullet} ) \big) \xrightarrow{\cong}
H^G_{\ast} \big( X^{\prime} ; C_{\ast}^{\HP} ( Z_{\bullet}^{\prime} ) \big), \\
H^G_{\ast} \big( X ; C_{\ast}^{\HN} ( Z_{\bullet} ) \big) \xrightarrow{\cong}
H^G_{\ast} \big( X^{\prime} ; C_{\ast}^{\HN} ( Z_{\bullet}^{\prime} ) \big)\,
\end{eqnarray*}
are isomorphisms.
\end{proposition}
\begin{proof}
There is a short exact sequence of chain complexes
\[
0 \to C_{\ast}^{\HH} (Z_{\bullet}) \to C_{\ast}^{\HC} (Z_{\bullet}) \to
C_{\ast}^{\HC} (Z_{\bullet}) [-2] \to 0\,,
\]
which is natural in $Z_{\bullet}$, see \cite[2.5.10 on pages 78--79]{Loday(1992)}. We use here the notation $C_{\ast}[r]$
for the chain complex which is shifted down $r$ steps, i.e.\ $(C_{\ast}[r])_n = C_{n+r}$. Since $\Tot^{\oplus}$ and
$\widetilde{C}_{\ast}^{\cell} ( X^?_+) \otimes_{\IZ \Or G} (-) $ are exact functors (we use here that $\widetilde
{C}_{\ast}^{\cell} ( X^?_+)$ is a free $\IZ \Or G$-module), the maps induced by $X \to X^{\prime}$ and $Z_\bullet
\to Z_\bullet^{\prime}$ lead to a short exact ladder diagram of chain complexes. The corresponding long
exact ladder in homology, the fact that $H_{\ast}^G ( X ; C_{\ast}^{\HH}( Z_{\bullet} ))$ and $H_{\ast}^G ( X ;
C_{\ast}^{\HC}( Z_{\bullet} ))$ are concentrated in non-negative degrees
and an easy inductive argument based on the Five-Lemma finish the proof for cyclic homology.

\smallskip

In order to prove the statement for periodic cyclic homology, one uses that the periodic cyclic complex can be considered
as the inverse limit of the tower of cyclic complexes
\begin{eqnarray} \label{tower}
\ldots \to C_{\ast}^{\HC}( Z_{\bullet} ) [4] \to C_{\ast}^{\HC}( Z_{\bullet} ) [2] \to
C_{\ast}^{\HC}( Z_{\bullet} ) [0]\,.
\end{eqnarray}
For $n\geq 0$, we have the following natural maps\,:
\begin{eqnarray*}
H_{n}^{G} \big(X ; C_{\ast}^{\HP} ( Z_{\bullet} ) \big) & \cong &
H_n \Big( \Tot^{\oplus} \big( \widetilde{C}_{\ast}^{\cell} ( X^?_+) \otimes_{\IZ \Or G} \lim_r C_{\ast}^{\HC}
(Z_{\bullet}) [2r] \big) \Big) \\[.3em]
& \to &
H_n \Big( \lim_r \, \Tot^{\oplus} \big( \widetilde{C}_{\ast}^{\cell} ( X^?_+) \otimes_{\IZ \Or G} C_{\ast}^{\HC}
(Z_{\bullet}) [2r] \big) \Big) \\[.3em]
& \to & \lim_r \,
H_n \Big( \Tot^{\oplus} \big( \widetilde{C}_{\ast}^{\cell} ( X^?_+) \otimes_{\IZ \Or G} C_{\ast}^{\HC} (Z_{\bullet})
[2r] \big) \Big) \\[.3em]
& \cong & \lim_r \, H_{n + 2r}^G \big( X ; C_{\ast}^{\HP} ( Z_{\bullet} ) \big)\,.
\end{eqnarray*}
In Lemma~\ref{lim-lemma} below, we show that the first map is an isomorphism if $X$ is a finite $G$-$CW$-complex.
The second map sits in a short exact $\lim^1$-$\lim$-sequence, because the maps in the tower \eqref{tower}
above are all surjective and the functors
$\widetilde{C}_{\ast}^{\cell}(X^?_+) \otimes_{\IZ \Or G} (-)$ and $\Tot^{\oplus}$ preserve surjectivity, compare
\cite[Theorem~3.5.8 on page~83]{Weibel(1994)}.
Since we already know the comparison result for the $\lim$- and $\lim^1$-terms involving cyclic homology, a Five-Lemma argument
yields the result for periodic cyclic homology.

\smallskip

It remains to prove the statement about negative cyclic homology.
There is a natural exact sequence of chain complexes \cite[5.1.4 on page~160]{Loday(1992)}
\[
0 \to
C_{\ast}^{\HN} ( Z_{\bullet} ) \to
C_{\ast}^{\HP} ( Z_{\bullet} ) \to
C_{\ast}^{\HC} ( Z_{\bullet} )[-2] \to
0\,.
\]
Again, one uses that $\widetilde{C}_{\ast}^{\cell}(X^?_+) \otimes_{\IZ \Or G} (-)$ and $\Tot^{\oplus}$ are exact
functors to produce a long exact ladder in homology and uses the Five-Lemma.
\end{proof}

In the previous proof we used the following statement.

\begin{lemma} \label{lim-lemma}
Suppose that $X$ is a finite $G$-$CW$-complex.
Then the natural map
{\small\begin{eqnarray*}
\;\;\Tot^{\oplus} \Big( \widetilde{C}_{\ast}^{\cell}(X^?_+) \otimes_{\IZ \Or G} \lim_r C_{\ast}^{\HC} (Z_{\bullet}) [2r] \Big)
\xrightarrow{\cong}\;
\lim_r\, \Tot^{\oplus} \big( \widetilde{C}_{\ast}^{\cell}(X^?_+) \otimes_{\IZ \Or G} C_{\ast}^{\HC} (Z_{\bullet}) [2r] \big)
\end{eqnarray*}}
\hspace*{-.3em}is an isomorphism.
\end{lemma}

\begin{proof}
There exists an exact sequence (by explicit construction of an inverse limit)
\[
0 \to
\lim_r C_{\ast}^{\HC}( Z_{\bullet} ) [2r]
\to
\prod_{r = 0}^{\infty} C_{\ast}^{\HC}( Z_{\bullet} ) [2r]
\to
\prod_{r = 0}^{\infty} C_{\ast}^{\HC}( Z_{\bullet} ) [2r]\,.
\]
As $\widetilde{C}_{\ast}^{\cell}(X^?_+) \otimes_{\IZ \Or G} (-)$ and $\Tot^{\oplus}$ are exact
functors, we see that it suffices to study the natural map
\[
\xymatrix{
\Tot^{\oplus} \big( \widetilde{C}_{\ast}^{\cell}(X^?_+) \otimes_{\IZ \Or G}
\prod_{r = 0}^{\infty} C_{\ast}^{\HC}( Z_{\bullet} ) [2r] \big)
\ar[d] \\
\prod_{r = 0}^{\infty} \Tot^{\oplus} \big( \widetilde{C}_{\ast}^{\cell}(X^?_+) \otimes_{\IZ \Or G}
 C_{\ast}^{\HC}( Z_{\bullet} ) [2r] \big)
}
\]
Let $C_\ast^{{\scriptscriptstyle\leq p}} \subseteq \widetilde{C}_{\ast}^{\cell}(X^?_+)$ be the $\Or G$-sub-complex
which agrees with $\widetilde{C}_{\ast}^{\cell}(X^?_+)$ up to dimension $p$ and is trivial in dimension $>p$. This
yields a finite filtration by our assumption on $X$. There is an induced map of filtered chain complexes
\[
\xymatrix{
F_\ast^p = \Tot^{\oplus} \big( C_\ast^{{\scriptscriptstyle\leq p}} \otimes_{\IZ \Or G}
\prod_{r = 0}^{\infty} C_{\ast}^{\HC}( Z_{\bullet} ) [2r] \big)
\ar@<1em>[d] \\
{}^{\prime}\!F_\ast^p = \prod_{r = 0}^{\infty} \Tot^{\oplus} \big( C_\ast^{{\scriptscriptstyle\leq p}} \otimes_{\IZ \Or G}
 C_{\ast}^{\HC}( Z_{\bullet} ) [2r] \big)
}
\]
and the induced chain map of filtration quotients $F_\ast^p / F_\ast^{p-1} \to {}^{\prime}\!F_\ast^p \big/
{}^{\prime}\!F_\ast^{p-1}$ can be identified with the composition
\[
\xymatrix{
\Tot^{\oplus} \big( \widetilde{C}_{p}^{\cell}(X^?_+) \otimes_{\IZ \Or G}
\prod_{r = 0}^{\infty} C_{\ast}^{\HC}( Z_{\bullet} ) [2r] \big)
\ar[d] \\
\Tot^{\oplus} \big( \prod_{r = 0}^{\infty} \widetilde{C}_{p}^{\cell}(X^?_+) \otimes_{\IZ \Or G}
C_{\ast}^{\HC}( Z_{\bullet} ) [2r] \big)
\ar[d] \\
\prod_{r = 0}^{\infty} \Tot^{\oplus} \big( \widetilde{C}_{p}^{\cell}(X^?_+) \otimes_{\IZ \Or G}
 C_{\ast}^{\HC}( Z_{\bullet} ) [2r] \big)
}
\]
because the tensor product with a fixed module over the orbit category, $\Tot^{\oplus}$ and $\prod_{r=0}^{\infty}$ all ``behave
well'' (in an obvious sense) with respect to taking quotients. The second map in the composition above is clearly an isomorphism.
The first map is an isomorphism because the assumption on $X$ implies that each $\widetilde{C}_{p}^{\cell}(X^?_+)$ is
a finitely generated free $\IZ \Or G$-module, compare \cite[page~167]{Lueck(1989)}. Since the filtrations are finite, this
concludes the proof.
\end{proof}

\begin{remark}
If we would only assume that $X$ is a $G$-$CW$-complex of finite type instead of being finite, then one would have the same
conclusion that the induced map of filtration quotients is an isomorphism for each $p$, as in the proof of Lemma~\ref{lim-lemma},
but the second filtration would not necessarily be exhaustive. The $0$th module of the complex
\[
\prod_{r = 0}^{\infty} \Tot^{\oplus} \big( \widetilde{C}_{\ast}^{\cell}(X^?_+) \otimes_{\IZ \Or G}
 C_{\ast}^{\HC}( Z_{\bullet} ) [2r] \big)
\]
would for instance contain the infinite product $\prod_{r=0}^{\infty} \widetilde{C}_{2r}^{\cell}(X^?_+)
\otimes_{\IZ \Or G} C_0^{\HC} ( Z_{\bullet} )$, whereas an element that is contained in ${}^{\prime}\!F_\ast^p$
for some $p$ has to be contained in the corresponding infinite direct sum.
\end{remark}

\bibliographystyle{plain}
\bibliography{dbdef,dbpub,dbpre,dbkhhextra}

\end{document}